\documentclass[aos,preprint]{imsart}
\RequirePackage[OT1]{fontenc}
\usepackage{amsthm,amsmath,amssymb,mathrsfs,graphicx}
\usepackage{verbatim,multirow,subfigure}
\usepackage{enumerate}
\usepackage{array}
\RequirePackage[numbers]{natbib}

\startlocaldefs

\theoremstyle{plain}

\newtheorem{theorem}{Theorem}[section]
\newtheorem{lemma}{Lemma}[section]
\newtheorem{proposition}{Proposition}[section]
\newtheorem{property}{Property}[section]
\newtheorem{corollary}{Corollary}[section]
\newcommand{\msum}{\mathop{\sum}}
\newcommand{\dmsum}{\mathop{\sum \sum}}
\newcommand{\xbar}{\overline{X}_n}
\newcommand{\twofam}[2]{#1^{(#2)}}
\newcommand{\mtxt}[1]{\mathcal{#1}_n}
\newcommand{\sammean}{\overline{X}_n'\overline{X}_n}
\newcommand{\varmn}{Var \left(\mtxt{M}\right)}
\newcommand{\trab}[2]{\widehat{ tr \left( \Gamma(#1) \Gamma(#2) \right)}}
\newcommand{\otrab}[2]{ tr \left( \Gamma(#1) \Gamma(#2) \right)}
\endlocaldefs
\setattribute{journal}{name}{}

\begin{document}

\begin{frontmatter}

% "Title of the paper"
\title{Mean vector testing for high dimensional dependent observations}
\runtitle{High dimensional mean testing}

\begin{aug}
\author{\fnms{Deepak} \snm{N. Ayyala}\corref{ayyala.1@osu.edu}\thanksref{m1,m2}\ead[label=e1]{ayyala.1@osu.edu}},
\author{\fnms{Junyong} \snm{Park}\thanksref{m2}\ead[label=e2]{junpark@umbc.edu}}
\and
\author{\fnms{Anindya} \snm{Roy}\thanksref{m2}
\ead[label=e3]{anindya@umbc.edu}}
%\ead[label=u1,url]{http://www.foo.com}}

%\thankstext{t1}{Some comment}
%\thankstext{t2}{First supporter of the project}
%\thankstext{t3}{Second supporter of the project}
\runauthor{Ayyala, Park and Roy}

\affiliation{The Ohio State University\thanksmark{m1} and University of Maryland Baltimore County\thanksmark{m2}}

\address{Address of the First author\\
 Department of Statistics\\
 The Ohio State University\\
 Columbus OH 43210\\
\printead{e1}}

\address{Address of the Second and Third authors\\
Department of Mathematics and Statistics\\\
University of Maryland Baltimore County\\
Baltimore MD 21250\\
\printead{e2}\\
\phantom{E-mail:\ }\printead*{e3}}

%\printead{u1}}
\end{aug}

\begin{abstract}
When testing for the mean vector in a high dimensional setting, it is generally assumed that the observations are independently and identically distributed. However if the data are dependent, the existing test procedures fail to preserve type I error at a given nominal significance level. We propose a new test for the mean vector when the dimension increases linearly with sample size and the data is a realization of an $M$-dependent stationary process. The order $M$ is also allowed to increase with the sample size. Asymptotic normality of the test statistic is derived by extending the central limit theorem result for $M$-dependent processes using two dimensional triangular arrays. Finite sample simulation results indicate the cost of ignoring dependence amongst observations.

\end{abstract}

\begin{keyword}[class=MSC]
\kwd[Primary ]{62H15}
\kwd{62F30,62M10}
\kwd[; secondary ]{91B84}
\end{keyword}

\begin{keyword}
\kwd{High dimension}
\kwd{asymptotic normality}
\kwd{triangular array}
\kwd{dependent data}
\kwd{mean vector testing}
\end{keyword}

\end{frontmatter}

\section{Introduction}

The advent of sophisticated data collecting methods are allowing researchers to measure a large number of variables simultaneously. In several applications in biomedical fields such as functional magnetic resonance imaging, microarray gene expressions and gene sequencing, the number of variables can be as high as several thousands. The cost of data collection is 
nevertheless still considerably high, restricting such studies to have no more than a few hundred samples. Standard multivariate techniques are no longer applicable to such data sets. This has created a huge surge in development of techniques for analyzing high-dimensional data and solving the so-called \textit{large p, small n} problems. One such problem of particular interest is testing for equality of mean vectors for two groups of observations. 

Let $\{X_{t} \}_{1 \leq t \leq n}$ be $p$-dimensional vectors with mean $\mu$ and covariance matrix $\Sigma$ and the problem of interest is to test
\begin{eqnarray}
H_0: \mu = 0 \hspace{6mm} \mbox{versus} \hspace{6mm} H_A: \mu \neq 0.
\label{eqn:onefammodel}
\end{eqnarray}
For two sample case,  let  $\{\twofam{X}{k}_t  \}, k = 1, 2$  be two independent groups of $p$-dimensional observations with mean $\twofam{\mu}{k}$ and covariance matrix $\twofam{\Sigma}{k}$ respectively. The hypothesis of interest is equality of the two mean vectors,
\begin{eqnarray}
 H_0: \twofam{\mu}{1} = \twofam{\mu}{2} \hspace{6mm} \mbox{versus} \hspace{6mm} H_A: \twofam{\mu}{1} \neq \twofam{\mu}{2}.
\label{eqn:twofammodel}
\end{eqnarray}
The traditional Hotelling's $T^2$ test is applicable only when the sample size is larger than the dimension and when the observations are independently and identically distributed, coming from a normal distribution. When one or both of these assumptions on the data fail to hold, Hotelling's $T^2$ cannot be applied. This calls for novel approaches in testing by relaxing the assumptions. 

Under the assumption that the observations are \textit{i.i.d} coming from a normal distribution, \cite{Dempster} proposed a approximate $F$-test statistic when the dimension is greater than the sample size. Significant progress has been made in this area by several researchers who have proposed test statistics which were proven to be asymptotically normal. The main idea behind these test statistics is constructing a function of the data whose expected value is a function of the population mean, which is equal to zero under the null hypothesis. \cite{BaiSaranadasa, Dempster} used the Euclidean norm of the sample mean and \cite{Srivastava,SrivastavaDu} used a variance-weighted quadratic product of the sample average. The functions used by \cite{ChenQin,ParkAyyala} are similar to those mentioned above, using only cross-product terms in the quadratic product of sample means and avoiding inner product terms which require additional conditions on the data.

Test statistics proposed by \cite{BaiSaranadasa, SrivastavaDu} directly constrain the rate of increase of dimension with respect to sample size while \cite{ChenQin,ParkAyyala} relax this direct relationship by imposing conditions on the dependence structure. Distributional assumption on the data required by \cite{BaiSaranadasa} have been relaxed by \cite{ChenQin,ParkAyyala,Srivastava2013,SrivastavaDu,Srivastava2013} and replaced by conditions on finiteness of fourth order moments. However all of these test statistics require the observations to be independent. When the observations are dependent with an unknown autocovariance structure, the Euclidean norm of the sample mean is no longer unbiased for the norm of the population mean. Hence the aforementioned test statistics will fail to preserve Type I error at any given nominal significance level. 

Structured dependency has been incorporated in various studies, but the dependence have been generally assumed on the variates. \cite{Zhong2013} have developed a test statistic based on a higher criticism test when the variates are weakly dependent and the observations are independent. Also, \cite{Cai2013} developed a test statistic which performs very well under sparse alternatives under dependence among the  variates. However, none of the studies incorporate a dependence structure amongst the observations. Testing for the mean vector under a dependence structure for the observations has been explored with a different perspective by several researchers. For example, \cite{Horvath1} studied stationarity of the mean vector when the data is assumed to be a realization of an $m-dependent$ strictly stationary model. A method for comparing two functional regression models was developed by \cite{Horvath3} wherein the model considered is a generalization of the factor model. \cite{Horvath2} addressed the problem of estimation of the mean for functional time 
series data under a weak stationary dependence structure. They also formulated an asymptotically normal test statistic for testing the equality of functional means of two populations. While the researches mentioned above address dependence between observations, none of them address testing for the hypothesis in \eqref{eqn:onefammodel}.

In this paper, we propose a test statistic which is based on the Euclidean norm of the sample mean, similar to \cite{BaiSaranadasa}. The proposed test statistic is constructed under the assumption that the rate of increase of dimension with respect to the sample size is linear and the observations are dependent. 

The rest of the paper is organized as follows. The factor model assumed for the data and construction of the test statistic for the one sample case are presented in Section \eqref{sec:onesample}. Asymptotic properties of the test statistic and construction of a ratio consistent estimator of the denominator of the test statistic are provided in Section \eqref{sec:asymptotics}. Extension of the test statistic to the two sample case is provided in Section \eqref{sec:twosample}. A finite sample simulation study is performed to evaluate the performance of the proposed test statistic and the results are presented in Section \eqref{sec:simulations}. Theoretical details and proofs of the theorems are given in the 
\nameref{sec:appendix}.

\section{One Sample case}
\label{sec:onesample}

We first address the problem for the one sample case. Let $X_t, 1 \leq t \leq n$ be $p$-dimensional vectors which follow a  $M$-dependent strictly stationary Gaussian process with mean $\mu$ and autocovariance structure given by $\Gamma(h), 0 \leq h \leq M$. Infinite moving average can approximate a large class of time series models. Thus, the models obtained  by allowing the order of the moving average to increase with the sample size constitute a rich class of models.  Previous studies for independent observations have emphasized the assumptions on dimensionality and the structure of covariance or correlation matrix. These conditions are required for deriving the asymptotic distribution of the test statistic. For dimensionality, Bai and Saranadasa \cite{BaiSaranadasa} (henceforth called $T_{BS}$) $( \frac{p}{n} \rightarrow c \in (0, \infty) )$  and Srivastava and Du \cite{SrivastavaDu} (henceforth called $T_{SD}$) $\left( n = O \left( p^{\zeta} \right), \frac{1}{2} < \zeta \leq 1 \right)$ assume a direct relationship between the sample size and dimension, whereas the test statistics proposed by Chen and Qin \cite{ChenQin} ($T_{CQ}$) and Park and Ayyala \cite{ParkAyyala} ($T_{PA}$) were constructed without a direct condition on $n$ and $p$. For the case of $M$-dependent model, we assume that the rate of increase of dimension with respect to the sample size is linear and the level of dependency $M$ increases at a polynomial rate,
\begin{eqnarray}
p&=& O(n), \label{eqn:dim_sample}  \\
M &=& O(n^{1/8}).   \label{eqn:M}
\end{eqnarray}
 The assumption  \eqref{eqn:M} on the order is necessary to have sufficiently large number of observations for estimating the autocovariance matrix at all lags. In addition to \eqref{eqn:dim_sample} and \eqref{eqn:M}, we impose conditions on the structure of the autocovariance matrix to avoid strong dependency amongst the variates. Define $\Omega_n$, the covariance matrix of $\overline X_n$ multiplied by the sample size as
\begin{eqnarray}
\Omega_n  \equiv n \, Cov(\overline X_n)  = \msum_{h = -M}^M \left(1 - \frac{|h|}{n} \right) \Gamma(h) = F(0).
 \label{eqn:omega}
\end{eqnarray}
where $F(0)$ is the spectral matrix evaluated at the zero frequency. 
To avoid considering arbitrary dependence structure amongst the variates, we put restrictions on the variance of the sample mean. The restrictions put constraints on the degree of dependence among the variates in relation to the degree of dependence among the observations. Such restrictions can be put on the spectral matrix. However, since we are primarily interested in a time domain test, we put restrictions on the autocovariance structure by imposing sparsity conditions on $\Gamma(h)$ through $\Omega_n$. 

For the case of independent observations, previous investigations imposed restrictions on covariance matrices in different ways to avoid arbitrary covariance matrices, especially the cases of highly correlated variables.
For example,  $T_{BS}$ and $T_{CQ}$ use   $\lambda_{\max} = o(\sqrt{tr\Sigma^2})$
and $tr(\Sigma^4) = o(tr^2(\Sigma^2))$ respectively, where $\lambda_{\max}$ is the maximum eigenvalue of
$\Sigma$. For $T_{SD}$ and $T_{PA}$, the correlation structure $\mathcal{R}$ is restricted by imposing $0 < \lim_{p \rightarrow \infty} \frac{tr \mathcal{R}^i}{p} < \infty, i = 1, 2, 3, 4$ and $tr (\mathcal{R}^4) = o( tr^2 \mathcal{R}^2)$ respectively.
For the case of dependent observational vectors, we put restrictions on the entire autocovariance structure to avoid strong dependency amongst the samples at all lags. To impose sparsity, it is assumed that for any set of four indices $-M \leq a, b, c, d \leq M$,
\begin{eqnarray}
 tr \left( \Gamma(a) \Gamma(b) \Gamma(c) \Gamma(d) \right) = o \left( tr^2 \left( \Omega_n^2 \right) \right),
 \label{eqn:condition1}
\end{eqnarray}
where the rate of decay is uniform for all $a,b,c$ and $d$. Together with \eqref{eqn:M}, the condition in \eqref{eqn:condition1} implies $tr ( \Omega_n^4 ) = \mathop{\sum }_{-M \leq  a,b,c,d  \leq M} w_{abcd} \, tr ( \Gamma(a) \Gamma(b) \Gamma(c) \Gamma(d))$ $=  o (n \,  tr^2 \left(  \Omega_n^2 \right) )$, where $w_{abcd} = \left( 1 - \frac{|a|}{n} \right)\left( 1 - \frac{|b|}{n} \right)\left( 1 - \frac{|c|}{n} \right)\left( 1 - \frac{|d|}{n} \right)$. This latter rate of decay assumption is similar to the sparsity conditions imposed by \cite{ChenQin,ParkAyyala,SrivastavaDu}.  However the assumption 
made in \eqref{eqn:condition1} is stronger as it requires uniform convergence over all combinations of indices. The more stringent condition is needed for addressing the dependence among the observations. 

We now present our proposed test statistic. For independent observational vectors, Bai and Saranadasa \cite{BaiSaranadasa} proposed
\begin{eqnarray}
 T_{BS} = \frac{n \overline{X}_n' \overline{X}_n - tr \left(S \right)}{\sqrt{ \frac{2 n(n-1)}{(n-2)(n+1)} \left( tr \left(S^2 \right) - \frac{ tr^2(S)}{n-1} \right) }}
\end{eqnarray}
where $\overline{X}_n$ is the sample average and $S = \frac{1}{n-1} \sum_i (X_i - \overline{X}_n) (X_i - \overline{X}_n)^{\prime}$ is the sample covariance matrix. They proved that $T_{BS}$ is asymptotically normal when $p/n \rightarrow c>0$. The main idea of $T_{BS}$ was to use $
E \left( n\overline{X}_n' \overline{X} - tr \left(S \right) \right) = \mu'\mu \geq 0$, which is equal to zero under the null hypothesis in \eqref{eqn:onefammodel} and non-zero under the alternative. However when the samples are dependent, the quantity $n\overline{X}_n' \overline{X} - tr \left(S \right)$ is no longer unbiased for $\mu' \mu$.
Since the autocovariance between observations causes a bias,  direct use of $T_{BS}$
is expected to fail in controlling a given size under $H_0$.
 A similar bias has been observed in all the previous test statistics such as $T_{CQ}$ for the moving average model. It is therefore not appropriate to use tests designed for independent observations when the samples are dependent.

The main idea of our test statistic is based on the Euclidean norm of the sample mean, $\sammean$, by correcting the bias due to the autocovariance structure. Since the observations are normally distributed, we have
\begin{eqnarray*}
 E(\sammean) &=& \frac{1}{n} \left[ \mu' \mu + \frac{1}{n} \msum_{h=-M}^M \msum_{k = -M}^M   tr\left(\Gamma(h-k) \right) \right]  =  \frac{1}{n} \mu' \mu + \frac{1}{n} tr(\Omega_n).
\end{eqnarray*}
If $\widehat{tr(\Omega_n)}$ is an unbiased estimator of $tr(\Omega_n)$, then the quantity $\mtxt{M} = \sammean - \frac{1}{n} \widehat{tr(\Omega_n)}$
will have zero mean under the null hypothesis. Denote by $\widehat{\Gamma}(h) = \frac{1}{n} \mathop{\sum}_{t = 1}^{n-h} (X_t - \xbar)(X_{t+h} - \xbar)'$ the biased sample estimator of the autocovariance matrix at lag $h$. When the dimension increases linearly with respect to sample size, this sample autocovariance matrix is no longer an asymptotically unbiased estimator for $\Gamma(h)$. Hence a plug-in estimator, $tr(\widehat{\Omega_n}) = \msum_{h = -M}^M \left(1 - \frac{|h|}{n} \right) tr \left(\widehat{\Gamma}(h) \right)$, will not work due to this bias.

We construct an unbiased estimator of $tr(\Omega_n)$ as follows, so that we have $E(\mtxt{M}) = \frac{1}{n} \mu' \mu$.
Let $\gamma$ be the $(M+1)\times 1$ vector consisting of $tr(\Gamma(h))_{0 \leq h \leq M}$ and $\widehat{\gamma}_n$ be the vector consisting of $tr(\widehat{\Gamma}(h))_{0 \leq h \leq M}$. The expected value of the vector consisting of these sample autocovariances can be calculated as $E \left( \widehat{\gamma}_n \right) = \Theta_n \gamma_n$,
where $\Theta_n$ is a $(M+1) \times (M+1)$ full-rank matrix whose elements are given by
\begin{eqnarray}
 \Theta_n(i,j) & = & \left(1 - \frac{i-1}{n} \right) \mathcal{J}(i,j) +  \left(1 - \frac{i-1}{n} \right)\left(1 - \frac{j-1}{n} \right) \frac{\left(2 - \mathcal{J}(i,1) \right)}{n}  \nonumber \\
           && - \frac{1}{n^2} \msum_{t=1}^{n-i+1} \msum_{s=1}^n \left( \mathcal{J}(|t-s|+1,j) + \mathcal{J}(|t+i-s-1|,j) \right),
 \label{eqn:thetas}
\end{eqnarray}
and $\mathcal{J}(a,b)$ is the indicator function, equal to one if $a$ and $b$ are equal and zero otherwise. Since $\Gamma(-h) = \Gamma(h)'$ for all lags, we can express $tr \left(\Omega_n \right)$ in terms of the autocovariance matrices at positive lags as $tr \left(\Omega_n \right) = b' \gamma$ where $b(1) = \frac{1}{n}$ and $b(i) = \frac{2}{n}\left(1 - \frac{i-1}{n} \right)$ for $2 \leq i \leq M$. An exactly unbiased estimator for $tr(\Omega_n)$ can be constructed as $  \widehat{tr(\Omega_n)} = \beta' \widehat{\gamma}_n =  b' \Theta_n^{-1} \widehat{\gamma}_n$. Using this unbiased estimator, it follows that the quantity $\mtxt{M} = \sammean - \frac{1}{n} \widehat{tr(\Omega_n)}$ will have expected value equal to $\frac{1}{n} \mu' \mu$. 

We propose $\mtxt{M}$ as the numerator and the test statistic is of the form $T_{new} = \frac{\mtxt{M}}{\sqrt{\widehat{\varmn}}}$, where $\widehat{\varmn}$ is a ratio-consistent estimator for the variance of the numerator. Expressing $tr(\widehat{\Gamma}(h))$ as a linear combination of inner product of $X_t$'s, we have
\begin{eqnarray}
\mathcal{M}_n & =&  \sammean - \frac{1}{n} \widehat{tr(\Omega_n)}= \msum_{t=1}^n \msum_{s=1}^n \pi_n(t,s) X_t' X_s, \hspace{4mm}
\label{eqn:mn}
\end{eqnarray}
where the coefficients can be calculated from \eqref{eqn:thetas} and the expression for $\widehat{tr \left(\Omega_n \right)}$  as
\begin{eqnarray*}
 \pi_n(t,s) & = & \frac{1}{n^2} \left( 1 - \frac{1}{n} \msum_{h=0}^M (n-h) \beta_n(h) \right) - \msum_{h=0}^M \Bigg[\frac{\beta_n(h)}{n} \mathcal{J}(i-j,h) \nonumber \\
 && - \frac{\beta_n(h)}{n^2} \left( \msum_{w=1}^{n-h} \mathcal{J}(i,w) + \msum_{w=h+1}^{n} \mathcal{J}(i,w) \right) \Bigg].
\end{eqnarray*}
Under this representation, the variance of $\mathcal{M}_n$ under the null hypothesis can be evaluated as
\begin{eqnarray}
 Var(\mathcal{M}_n) = 2 \msum_{h=-M}^M \msum_{k=-M}^M \xi_n(h,k) tr \left(\Gamma(h) \Gamma(k) \right) = \xi_n' \lambda,
\label{eqn:VarM}
\end{eqnarray}
where $\xi_n$ and $\lambda$ are $(M+1)^2 \times 1$ vectors, which are vector forms of $\left(M+1 \right) \times \left( M+1 \right)$ matrices whose elements are given by $\lambda(h,k) = tr \left(\Gamma(h) \Gamma(k) \right)$ and $\xi_n(h,k) = \msum_{i,j=1}^n \pi_n(i,j) \pi_n(j-h,i+k)$ respectively. Length of the vectors $\xi_n$ and $\lambda$ is taken as $(M+1)^2$ instead of $(2M+1)^2$ because of the properties of the autocovariance matrices, $tr \left(\Gamma(a) \Gamma(b) \right) = tr \left(\Gamma(b) \Gamma(a) \right) = tr \left(\Gamma(-a) \Gamma(-b) \right)$ $= tr \left(\Gamma(-b) \Gamma(-a) \right)$. Therefore, the number of terms to be estimated in $\lambda$ is reduced by a factor of $\frac{1}{4}$ and the terms in $\xi_n$ are adjusted accordingly for multiplicity. 

\section{Asymptotic properties of the test statistic}
\label{sec:asymptotics}

In this section, we derive the asymptotic null distribution of $\mtxt{T} = \frac{\mtxt{M}}{\sqrt{\varmn}}$ and construct a ratio consistent estimator for $\varmn$. 

\begin{theorem}
\label{thm:clt}
 When the conditions specified in \eqref{eqn:dim_sample}, \eqref{eqn:M} and  \eqref{eqn:condition1} are satisfied, then under the null hypothesis
\begin{eqnarray*}
 \mathcal{T}_n = \frac{\mtxt{M}}{\sqrt{Var\left(\mtxt{M}\right)}} \xrightarrow{\mathcal{D}} \mathcal{N}(0,1) \mbox{ as } n \rightarrow \infty
\end{eqnarray*}
where $\xrightarrow{\mathcal{D}}$ means convergence in distribution.
\end{theorem}
\noindent \textbf{Proof}: See \nameref{sec:appendix}.  \qed

The main idea of the proof is constructing a $n \times n$ matrix $\mathcal{A}$ whose elements are the summands of $\mtxt{M}$ from \eqref{eqn:mn}. For univariate $M$-dependent processes, \cite{Berk,Diananda,Hoeffding,RomanoWolf} have derived central limit theorems using triangular arrays. Novelty of our proof lies in the extension of triangular array argument to two dimensional arrays. 

Let $\widehat{\varmn}$ be a ratio consistent estimator for $\varmn$, that is $\frac{\widehat{Var \left(\mtxt{M} \right)}} {Var(\mathcal{M}_n)} \stackrel{p}{\rightarrow} 1$, where $\stackrel{p}{\rightarrow} 1$ means convergence in probability. Before constructing a ratio-consistent estimator for the variance, an equivalent representation of $\varmn$ and its limiting expression are described in the following proposition. 
\begin{proposition}
The coefficients $\xi_n(h,k)$'s in \eqref{eqn:VarM} can be expressed as
 \begin{align*}
  \xi_n(h,k) = \frac{1}{n^2} \left(1 + \chi_n \right) \left(1 - \frac{|h|}{n} \right) \left(1 - \frac{|k|}{n} \right)(1 + o(1)), \,  -M \leq h,k \leq M,
 \end{align*}
 where $\chi_n = \frac{1}{n} \msum_{h=1}^M \left(1 - \frac{h}{n} \right)^2 = o(1)$. Therefore the variance can be expressed as $ Var(\mathcal{M}_n)  = \frac{2}{n^2} \left(1 + \chi_n \right) \, tr \left[ \msum_{h=-M}^M \left(1 - \frac{|h|}{n} \right)\Gamma(h) \right]^2 \left(1 + o(1) \right)$ $= \frac{2}{n^{2}} tr \left(\Omega_n^2 \right)  + o \left(n^{-2} tr \left(\Omega_n^2 \right) \right)$.
\label{prop:simplexhis}
\end{proposition}
{\bf Proof :} See \nameref{sec:appendix}.  \qed

We now propose our estimator for the variance. To begin with, a plug-in estimator has two major drawbacks. Firstly, it is not unbiased for the variance because for any $-M \leq a,b \leq M$, the expected trace of product of sample autocovariance matrices at lags $a$ and $b$ is a function of the entire autocovariance structure, $ E \left[ tr \left( \widehat{\Gamma(a)} \widehat{\Gamma(b)} \right) \right] = \dmsum_{i,j = -M}^M \tau_n(i,j) tr \left( \Gamma(i) \Gamma(j) \right)$,
where the coefficients $\tau_n(i,j)$ can be calculated by simplifying the expression for the expected value. An unbiased plug-in estimator can therefore be constructed as $\widehat{\varmn} = \xi' \tau_n^{-1} \widehat{\lambda}$, where the $(M+1)^2 \times (M+1)^2$ matrix $\tau_n$ is of full rank. However, establishing ratio-consistency of this unbiased estimator is arduous. To be able to  demonstrate that the variance of the estimator constructed converges to zero, that is $Var \left(\widehat{\varmn} \right) = o \left(\varmn^2 \right)$, additional conditions are required on the autocovariance structure. Secondly, the above estimator is not guaranteed to be positive since it is not based on the quadratic form of the vector of traces.

We construct the estimator $\widehat{\varmn}$ as follows. For any $-M \leq a,b \leq M$, we first define an estimator of the trace of product of autocovariance matrices as
\begin{align}
 \widehat{tr \left(\Gamma(a) \Gamma(b) \right)} = \frac{1}{n_{a,b}} \mathop{\sum \sum}_{t,s \in \mathcal{A}(a,b)} \left(X_{t+a} - \overline{X}_n^* \right)^{\prime} X_{s} \left(X_{s+b} - \overline{X}_n^* \right)^{\prime} X_{t},
  \label{eqn:traceest}
\end{align}
where $\mathcal{A}(a,b) = \left\{ (t,s) : |t+a-s-b| > M, |t-s| > M \right\}$, $n_{a,b}$ is the cardinality of $\mathcal{A}(a,b)$ and $\overline{X}_n^*$ is the average of all observations in $\left\{ X_w \right\}_{w \in \mathcal{B}(t,s)}$  for any $(t,s) \in \mathcal{A}(a,b)$ and $\mathcal{B}(t,s) = \{ 1 \leq i \leq n : \min\{|i-t|, |i-s|, |i-t-a|, |i-s-b|\} > M \}$. For notational convenience, unless otherwise stated, $\mathcal{A}(a,b)$ will be denoted by $\mathcal{A}$ and $\mathcal{B}(t,s)$ will be denoted by $\mathcal{B}$ henceforth.

The above method of estimating $tr  \left(\Gamma(a) \Gamma(b) \right)$ is an extension of the idea proposed by Chen and Qin \cite{ChenQin} in that
 $\mathcal{A}$ and $\mathcal{B}$ are constructed noting that two observations $X_t$ and $X_s$ are independent when they are at least $M$ time points apart. The above construction of $\mathcal{A}$ and $\mathcal{B}$ ensures that the leading term of the proposed estimator is consistent for $\varmn$. The estimator constructed using $\widehat{tr \left(\Gamma(a) \Gamma(b) \right)}$ is also required to be ratio consistent for non-zero $\mu$ which belong to a local alternative. Note that only two of the terms in the product of quadratics in \eqref{eqn:traceest} are centered to avoid dealing with the term $ \left(\overline{X}_n^{* '} \overline{X}_n^*\right)^2$, proving the consistency of which requires additional conditions on $\Omega_n$.

For any $-M \leq a, b \leq M$, since the autocovariance matrices at non-zero lags are not necessarily positive definite and symmetric, it is not guaranteed that $tr\left( \Gamma(a) \Gamma(b) \right)$ is non-zero even though both $\Gamma(a)$ and $\Gamma(b)$ non-zero. Therefore, it would be inexact to establish \eqref{eqn:traceest} as a ratio consistent estimator for $tr\left( \Gamma(a) \Gamma(b) \right)$ for all $-M \leq a,b \leq M$. However, $\varmn$ is a nonzero quantity and we can construct a ratio consistent estimator for $\varmn$. Further to allow for ratio-consistency of the variance estimator for non-zero mean, we restrict the mean to a local alternative given by 
\begin{eqnarray}
 \label{eqn:alternativemean}
 \mu' \left[\Gamma(a) \Gamma(-a) \right]^{\frac{1}{2}} \mu = o \left( n^{-1} tr \left(\Omega_n^2\right) \right) \mbox{   for   } -M \leq a \leq M.
\end{eqnarray}

To better appreciate this condition, set the lag equal to zero to obtain $\mu' \Gamma(0) \mu = o \left( n^{-1} tr \left(\Omega_n^2\right) \right)$, which is similar to the local alternative condition specified in literature. To summarize, the following theorem provides a ratio consistent estimator of $\varmn$.
\begin{theorem}
 \label{thm:variance}
 If the conditions in Theorem \ref{thm:clt} and \eqref{eqn:alternativemean} hold and the estimator of the trace $ \widehat{tr \left(\Gamma(a) \Gamma(b) \right)}$ is as mentioned in \eqref{eqn:traceest}, then
 \begin{eqnarray*}
  \widehat{Var \left(\mtxt{M} \right)} = \mathop{\sum\sum}_{a,b = -M}^{M} \xi_n(a,b) \widehat{tr \left(\Gamma(a) \Gamma(b) \right)},
 \end{eqnarray*}
is ratio consistent for $Var \left( \mtxt{M} \right)$, i.e.  $\frac{\widehat{Var \left(\mtxt{M} \right)} }{Var \left( \mtxt{M} \right)} \stackrel{p}{\rightarrow} 1$.
\end{theorem}
{\bf Proof :} See \nameref{sec:appendix}.   \qed

Combining Theorems \ref{thm:clt} and \ref{thm:variance}, we have asymptotic normality of our proposed test.  We state this result in the following corollary.

\begin{corollary} Under the assumptions specified in Theorems \eqref{thm:clt} and \eqref{thm:variance} the test statistic is asymptotically normal under the null hypothesis \eqref{eqn:twofammodel},
\begin{eqnarray*}
T_{new} = \frac{\mtxt{M}}{\sqrt{\widehat{Var \left(\mtxt{M} \right)}}} \xrightarrow{\mathcal{D}} \mathcal{N}(0,1) \, \, as \, \, n \rightarrow \infty.
\end{eqnarray*}
\end{corollary}
\noindent {\bf Proof :} The convergence in distribution of $T_{new}$ is straightforward from
\begin{eqnarray*}
 \frac{\mtxt{M}}{\sqrt{\widehat{Var \left(\mtxt{M} \right)}}}  =  \frac{\sqrt{{Var \left(\mtxt{M} \right)}}}   {\sqrt{\widehat{Var \left(\mtxt{M} \right)}}}  \frac{\mtxt{M}}{\sqrt{{Var \left(\mtxt{M} \right)}}}  \xrightarrow{\mathcal{D}}  \mathcal{N}(0,1)
\end{eqnarray*}
due to Theorem \ref{thm:clt}, Theorem \ref{thm:variance} and Slutsky's Theorem.
\qed

\bigskip

\noindent \textbf{Rejection rule and asymptotic power}

Though the alternative in \eqref{eqn:onefammodel} is two sided, the expected value of the numerator in $T_{new}$ under the alternative is equal to $ \mu' \mu$ which is always positive. Hence the hypothesis can be equivalently represented as
\begin{eqnarray*}
 H_0: \mu' \mu = 0 \hspace{6mm} \mbox{versus} \hspace{6mm} H_A: \mu' \mu > 0,
\end{eqnarray*}
which indicates the alternative is in fact one sided. The rejection region is therefore one sided, based on the upper tail of the normal distribution as opposed to a two-tailed rejection region. This reasoning applies to all test statistics based on the Euclidean norm  and as explained in \cite{BaiSaranadasa,ChenQin,ParkAyyala,SrivastavaDu}, the null hypothesis is rejected at significance level $\alpha$ when $T_{new} >z_{\alpha}$, where $z_{\alpha}$ is the $(1-\alpha)$ percentile of the standard normal distribution.

Under the local alternative condition specified in \eqref{eqn:alternativemean}, the estimator of variance is consistent even for non-zero mean. Asymptotic normality of
$  \frac{\mtxt{M} - \mu'\mu}{\sqrt{\widehat{Var \left(\mtxt{M} \right)}}} \xrightarrow{\mathcal{D}} N(0,1)$ follows from $E(\mtxt{M})= \mu'\mu$.
From this result, the asymptotic power function $\beta_n \left(\mu \right)$ of the proposed test statistic under a local alternative can be expressed as
\begin{eqnarray}
 \beta_n \left(\mu \right) - \Phi \left( -z_{\alpha} + \frac{ n \mu' \mu}{\sqrt{2 tr \left(\Omega_n^2 \right)}} \right) \rightarrow 0 \, \, \mbox{  as  } \, \, n \rightarrow \infty,
 \label{eqn:asymppower}
\end{eqnarray}
where $\Phi$ is the cumulative distribution function of the standard normal distribution.

%---------------------------------------------------------Section-------------------------------------------------
\section{Two Sample Case}
\label{sec:twosample}
In this section, we extend the test statistic proposed for one sample case to the two sample case. Let $\left\{ \twofam{X}{1}(t) \right\}$ and $\left\{ \twofam{X}{2}(t) \right\}$ be two families of $p$-variate vectors coming from independent $M$-dependent Gaussian processes. The means and autocovariance structures of the two groups are given by $\twofam{\mu}{1}$ and $\twofam{\mu}{2}$ and $\{\twofam{\Gamma}{1}(h) \}_{0 \leq h \leq M}$ and $\{\twofam{\Gamma}{2}(h) \}_{0 \leq h \leq M}$ respectively. The autocovariance structures of the two groups need not necessarily be the same. Most of the studies on independent observations impose the same covariance structure on both samples. For example, \cite{BaiSaranadasa} assumes
 the homogeneous covariance matrix and  \cite{ParkAyyala} and \cite{Srivastava} assume the same correlation matrix for the two sample case while \cite{ChenQin} does not demand the homogeneous covariance matrix. The hypothesis of interest is equality of means of the two samples,
\begin{eqnarray}
 H_0: \twofam{\mu}{1} = \twofam{\mu}{2} \hspace{6mm} \mbox{versus} \hspace{6mm} H_A: \twofam{\mu}{1} \neq \twofam{\mu}{2}.
 \label{eqn:twofamhyp}
\end{eqnarray}

Let $\twofam{\Omega_{n_1}}{1}$ and $\twofam{\Omega_{n_2}}{2}$ be $p \times p$ matrices defined as in \eqref{eqn:omega}. Then the expected value of the Euclidean norm of difference of the sample means can be evaluated as
\begin{align*}
 E \left[ \left(\twofam{\overline{X}}{1}- \twofam{\overline{X}}{2} \right)^{\prime} \left(\twofam{\overline{X}}{1}- \twofam{\overline{X}}{2} \right) \right] &=& \left(\twofam{\mu}{1}- \twofam{\mu}{2} \right)^{\prime} \left(\twofam{\mu}{1}- \twofam{\mu}{2} \right) \\
 && + \qquad{} tr \left( \frac{1}{n_1} \twofam{\Omega_{n_1}}{1} + \frac{1}{n_2} \twofam{\Omega_{n_2}}{2} \right).
\end{align*}
Let $n = min(n_1, n_2)$. The test statistic proposed for the one sample case can be modified for the two sample case as $T_{new} = \frac{\mtxt{M}}{\sqrt{ \widehat{Var \left(\mtxt{M} \right)}}}$, where
\begin{eqnarray*}
 \mtxt{M} = \left(\twofam{\overline{X}}{1}- \twofam{\overline{X}}{2} \right)' \left(\twofam{\overline{X}}{1}- \twofam{\overline{X}}{2} \right) - \left( \frac{1}{n_1} \widehat{tr \left(\twofam{\Omega_{n_1}}{1} \right)} + \frac{1}{n_2} \widehat{tr \left(\twofam{\Omega_{n_2}}{2} \right)}\right),
\end{eqnarray*}
and the estimators $\widehat{tr \left(\twofam{\Omega_{n_1}}{1} \right)}$ and $\widehat{tr \left(\twofam{\Omega_{n_2}}{2} \right)}$ are exactly unbiased and as constructed for the one sample case.
To establish asymptotic normality and construct a ratio consistent estimator of the variance, the assumptions made in Section \ref{sec:onesample} need to be extended as follows. For the two sample case, sparsity condition on the autocovariance structures and the local alternative condition for the means are 
\begin{align}
 tr \left( \twofam{\Gamma}{w_1} (a) \twofam{\Gamma}{w_2} (b) \twofam{\Gamma}{w_3} (c) \twofam{\Gamma}{w_4} (d) \right) = o \left( tr^2 \left(\twofam{\Omega_{n_1}}{1} + \twofam{\Omega_{n_1}}{2} \right)^2 \right), \label{eqn:twofamass1}
 \end{align}
 \begin{align}
 (\twofam{\mu}{1}- \twofam{\mu}{2} )'[\twofam{\Gamma}{w}(a) \twofam{\Gamma}{w}(-a) ]^{\frac{1}{2}} (\twofam{\mu}{1}- \twofam{\mu}{2} ) = o (n^{-1} tr (\twofam{\Omega_{n_1}}{1} + \twofam{\Omega_{n_1}}{2} )^2 ), \label{eqn:twofamass2}
\end{align}
for all $\left(w_1,w_2,w_3,w_4 \right) \in \{1,2\}^4$ and $-M \leq a,b,c,d \leq M$. 

Express $\mtxt{M} = \twofam{\mtxt{M}}{1} + \twofam{\mtxt{M}}{2} - 2 \twofam{\mtxt{M}}{3}$, where the three terms are defined as $\twofam{\mtxt{M}}{1} = \overline{X}^{(1) \prime} \overline{X}^{(1)} - \frac{1}{n_1} \widehat{tr \left(\twofam{\Omega_{n_1}}{1} \right)}$, $\twofam{\mtxt{M}}{2} = \overline{X}^{(2) \prime} \overline{X}^{(2)} - \frac{1}{n_2} \widehat{tr \left(\twofam{\Omega_{n_2}}{2} \right)}$ and $\twofam{\mtxt{M}}{3} = \overline{X}^{(1) \prime} \overline{X}^{(2)}$ respectively. Using this form for $\mtxt{M}$ and under $M$-dependence Gaussian property of the model, it is straightforward to show that the quantities $\twofam{\mtxt{M}}{1}$, $\twofam{\mtxt{M}}{2}$ and $\twofam{\mtxt{M}}{3}$ are independent and the variance under the null hypothesis is given by
\begin{align*}
 Var \left(\mtxt{M} \right) &= Var \left(\twofam{\mtxt{M}}{1} \right) + Var \left(\twofam{\mtxt{M}}{2} \right) + 4 \, Var \left(\twofam{\mtxt{M}}{3} \right) \\
 &= 2 \dmsum_{h, k=-M}^M \bigg[ \xi_{n_1}(h,k) tr \left( \twofam{\Gamma}{1}(h) \twofam{\Gamma}{1}(k) \right)   + \xi_{n_2}(h,k) tr \left( \twofam{\Gamma}{2}(h) \twofam{\Gamma}{2}(k) \right) \bigg] \\
 & + \frac{4}{n_1 n_2} tr \left( \twofam{\Omega_{n_1}}{1} \twofam{\Omega_{n_2}}{2} \right).
\end{align*}

The first two terms in the above expression are equal to $\varmn$ as in the one sample. Therefore we can use Theorem \ref{thm:variance} to construct ratio consistent estimators for $Var \left(\twofam{\mtxt{M}}{1} \right)$ and  $Var \left(\twofam{\mtxt{M}}{2} \right)$. For the third term, we consider an estimator of the form
\begin{eqnarray*}
 \widehat{ tr \left( \twofam{\Omega_{n_1}}{1} \twofam{\Omega_{n_2}}{2} \right)} &= &\dmsum_{a,b=-M}^M \left(1 - \frac{|a|}{n_1} \right)\left(1 - \frac{|b|}{n_2} \right) \widehat{tr \left(\twofam{\Gamma}{1}(a) \twofam{\Gamma}{2}(b) \right)}. 
 \end{eqnarray*}
 As discussed in Section \ref{sec:onesample}, using a plug-in estimator with $\widehat{tr \left(\twofam{\Gamma}{1}(a) \twofam{\Gamma}{2}(b) \right)} = tr \left( \widehat{\twofam{\Gamma}{1}(a)} \widehat{\twofam{\Gamma}{2}(b)} \right)$ is not feasible. An estimator constructed similar to \eqref{eqn:traceest} instead of a plug-in estimator is used,
 \begin{align}
 \widehat{tr \left(\twofam{\Gamma}{1}(a) \twofam{\Gamma}{2}(b) \right)} = \frac{1}{n_{a,b}} \msum_{t=1}^{n_1-a} \msum_{s=1}^{n_2-b} (\twofam{X}{1}_{t+a} - \overline{\twofam{X}{1}}^* )^{\prime} \twofam{X}{2}_{s} (\twofam{X}{2}_{s+b} - \overline{\twofam{X}{2}}^* )^{\prime} \twofam{X}{1}_{t},
 \label{eqn:omega12}
\end{align}
where $\overline{\twofam{X}{w}}^*$ is the average of the set of observations with indexing set given by $\twofam{\mathcal{B}}{w}(t) =$ $ \left\{ 1 \leq i \leq n_w : |i-t|>M, |i-t-a|>M \right\}$ for $w=1,2$. As explained in the one sample case, the indexing sets are constructed so that $\overline{\twofam{X}{1}}^*$  and $\overline{\twofam{X}{2}}^*$are independent of the other terms in the expression. The following theorem provides the asymptotic distribution of $\mtxt{T}$ and a ratio consistent estimator for the variance.
\begin{theorem}
 \label{thm:twosample}
  If the conditions specified in \eqref{eqn:dim_sample}, \eqref{eqn:M}, \eqref{eqn:twofamass1} and \eqref{eqn:twofamass2}  are satisfied, then $T_{new}$ is asymptotically normal under the null hypothesis,
  \begin{eqnarray}
   T_{new} = \frac{\mtxt{M}}{\sqrt{\widehat{\varmn}}} \xrightarrow{\mathcal{D}} \mathcal{N} \left(0, 1\right),
  \end{eqnarray}
where the estimator of the variance is constructed using the estimators proposed in \eqref{eqn:traceest} and \eqref{eqn:omega12},
\begin{align*}
 \widehat{\varmn} &= 2 \dmsum_{h, k=-M}^M \bigg[ \xi_{n_1}(h,k) \widehat{tr \left( \twofam{\Gamma}{1}(h) \twofam{\Gamma}{1}(k) \right)} \\
 & +  \xi_{n_2}(h,k) \widehat{ tr \left( \twofam{\Gamma}{2}(h) \twofam{\Gamma}{2}(k) \right)} \bigg] + \frac{4}{n_1 n_2} \widehat{tr \left( \twofam{\Omega_{n_1}}{1} \twofam{\Omega_{n_2}}{2} \right)}.
\end{align*}

\end{theorem}
\noindent \textbf{Proof}: See \nameref{sec:appendix}.  \qed

Under the local alternative condition \eqref{eqn:twofamass2}, the asymptotic power of the two-sample test can be calculated as
\begin{eqnarray}
\;\;\;\;\;\;\;\;\;\;\; \beta_n \left(\twofam{\mu}{1} - \twofam{\mu}{2} \right) - \Phi \left(-z_{\alpha} + \frac{\left(\twofam{\mu}{1}- \twofam{\mu}{2} \right)' \left(\twofam{\mu}{1}- \twofam{\mu}{2} \right)}
 {\sqrt{2 tr \left(\frac{1}{n_1} \twofam{\Omega_n}{1} + \frac{1}{n_2} \twofam{\Omega_n}{2} \right)^2  }} \right) \rightarrow 0.
 \label{eqn:twosamasymppower}
\end{eqnarray}

%---------------------------------------------------------Section-------------------------------------------------
\section{Numerical Study}
\label{sec:simulations}
To evaluate the performance of the proposed test statistic against other test procedures mentioned, the size and power of the test statistics were compared for different levels of dependency.
We compare our proposed test with test procedures for independent observations such as $T_{CQ}$, $T_{BS}$, $T_{SD}$ and $T_{PA}$ and highlight on demonstrating the failure of these test statistics when observations are dependent while the proposed test obtains reasonable performance.

To generate samples from a $M$-dependent stationary Gaussian process, we considered an extension of the factor model used in \cite{ChenQin,ParkAyyala,SrivastavaDu}. The samples are assumed to be coming from a factor model with stationary mean
\begin{eqnarray}
 X_t = \mu + \mathop{\sum}_{h=0}^{M} A_h \varepsilon_{t-h},
 \label{eqn:model}
\end{eqnarray}
where $\left\{A_h\right\}_{0 \leq h \leq M}$ are unknown $p \times m$ matrices which determine the autocovariance structure.
To ensure the autocovariance matrix at lag zero has full rank, we assume $m > p$. The elements of the vector $\varepsilon_t = \left( \varepsilon_{ti} \right)_{1 \leq i \leq m}$ are assumed to be generated independently and identically from  a $m$-dimensional normal distribution with zero mean and covariance matrix equal to identity. Covariance matrix of the error vectors can be taken to be equal to the identity matrix without loss of generality since if the variance is equal to $\Sigma$, then the model will be equivalent to (\ref{eqn:model}) with $A_h = A_h \Sigma^{\frac{1}{2}}$. The model (\ref{eqn:model}) with the normality of $\varepsilon_t$ implies $X_t \sim N_p(\mu, \Gamma(0))$ and $Cov(X_t, X_{t+h}) = \Gamma(h)$  for $0 \leq h \leq M$, where $\Gamma(h) = \mathop{\sum}_{k=0}^{M-h} A_k  A_{k+h}', \, \, \, 0 \leq h \leq M$.

\subsection{One sample case}
\label{sec:oscase}
We consider the testing problem for one sample case under four different situations depending on different autocovariance matrices and different configurations of $(n,p,M)$.
The data are generated from the model in \eqref{eqn:model} where the mixing matrices and the error covariance matrix are given by
\begin{align}
 A_h(i,j) = \left\{
     \begin{array}{cc}
       \frac{\phi_1}{h |i-j|^2} & \mbox{if } |i-j| \leq p w, \\
       0 & \mbox{if } |i-j| > p w,
     \end{array}
   \right.
   &
   \Sigma(i,j) = \left\{
     \begin{array}{cc}
       \frac{\sqrt{\sigma_i \sigma_j} \phi_2}{|i-j|^2} & \mbox{if } |i-j| \leq p w, \\
       0 & \mbox{if } |i-j| > p w,
     \end{array}
   \right.
   \nonumber
\end{align}
respectively, where $0 \leq w\leq 1$ controls sparsity and the parameters $\phi_1, \phi_2$ and $\sigma$ control dependence
amongst the variates. We considered four models with different set of values for the parameters. For the four models named Model I through Model IV, the level of dependency $M$ is fixed to be 0, 1, 2 and 3 respectively. For Model I, the parameters were fixed as $\left(\frac{p}{n},\phi_1,\phi_2,w \right) = \left(4, 0.2, 0.3, 0.9\right)$ while the parameter values for Models II through IV are $\left(1, 0.6, 0.4, 0.8 \right), \left(2, 0.6, 0.6, 0.8 \right)$ and $\left( 3, 0.6, 0.3, 0.8 \right)$ respectively.

Further for models III and IV, the elements of the mixing matrix were as defined above with $A_h(i,j) = \frac{h \phi_1}{|i-j|^2}  \mbox{  for  } |i-j| \leq p$. For all the models, the power was studied when the elements of the mean vector were generated independently with $\mu_i \sim \frac{1}{\sqrt[4]{p}} U(2,3)$ or $\mu_i \sim \frac{1}{\sqrt{p}} U(2,3)$ and the sample size was taken to be equal to $40, 60, 80$ or $100$.

\begin{table}[htpb]
\renewcommand\arraystretch{1.4}
 \centering
 \begin{tabular}{|p{3.3em}|>{\centering\arraybackslash}p{1.em}|*{5}{>{\centering\arraybackslash}p{2.2em}}|*{5}{>{\centering\arraybackslash}p{2.2em}}|}
\multicolumn{2}{c}{}& \multicolumn{5}{c}{Model I} &  \multicolumn{5}{c}{Model II} \\
\hline
  & n & $T_{new}$ & $T_{BS}$ & $T_{CQ}$ & $T_{PA}$ & $T_{SD}$ & $T_{new}$ & $T_{BS}$ & $T_{CQ}$ & $T_{PA}$ & $T_{SD}$ \\
  \hline
  \multirow{4}{*}{Size}
& 40  &0.061&  	0.097&  	0.061&  	0.071&  0.065&  0.076&  0.442&  0.405&  0.431&  0.358\\
& 60  &0.063&  	0.088&  	0.063&  	0.064&  0.061&  0.070&  0.533&  0.507&  0.526&  0.454\\
& 80  &0.055&  	0.079&  	0.055&  	0.060&  0.055&  0.070&  0.611&  0.589&  0.602&  0.534\\
& 100 &0.054&  	0.074&  	0.054&  	0.055&  0.050& 	0.068&  0.674&  0.658&  0.670&  0.605\\
\hline
  \multirow{4}{*}{Power I}
& 40  &0.989&  	0.994&  0.989&  0.991&  0.989&  0.242&  0.678&  0.649&  0.672&  0.620\\
& 60  &0.999&  	0.999&  0.999&  0.999&  0.999&  0.282&  0.786&  0.767&  0.787&  0.751\\
& 80  &1&  	1&  	1&  	1&  	1&  	0.319&  0.858&  0.846&  0.86&  	0.821\\
& 100 &1&  	1&  	1&  	1&  	1&  	0.35&  	0.906&  0.897&  0.908&  0.881\\
\hline
  \multirow{4}{*}{Power II}
& 40  &1&  	1&  	1&  	1&  	1&  	0.818&  	0.934&  	0.929&  	0.933&  	0.99\\
& 60  &1&  	1&  	1&  	1&  	1&  	0.989&  	1&  	1&  	1&  	1\\
& 80  &1&  	1&  	1&  	1&  	1&  	0.999&  	1&  	1&  	1&  	1\\
& 100 &1&  	1&  	1&  	1&  	1&  	1&  	1&  	1&  	1&  	1\\
\hline
\multicolumn{2}{c}{}& \multicolumn{5}{c}{Model III} &  \multicolumn{5}{c}{Model IV} \\
\hline
  & n & $T_{new}$ & $T_{BS}$ & $T_{CQ}$ & $T_{PA}$ & $T_{SD}$ & $T_{new}$ & $T_{BS}$ & $T_{CQ}$ & $T_{PA}$ & $T_{SD}$ \\
  \hline
  \multirow{4}{*}{Size}
& 40  &0.072&  	0.929&  	0.911&  	0.937&  	0.902&  	0.060&  	0.998&  	0.995&  	0.998&  	0.995\\
& 60  &0.071&  	0.979&  	0.974&  	0.983&  	0.969&  	0.063&  	1&  	1&  	1&  	1\\
& 80  &0.068&  	0.996&  	0.995&  	0.996&  	0.993&  	0.062&  	1&  	1&  	1&  	1\\
& 100 &0.065&  	0.999&  	0.999&  	0.999&  	0.998&  	0.058&  	1&  	1&  	1&  	1\\
\hline
  \multirow{4}{*}{Power I}
& 40  &0.125&  	0.952&  	0.937&  	0.959&  	0.916&  	0.084&  	0.997&  	0.997&  	0.998&  	0.996\\
& 60  &0.135&  	0.989&  	0.985&  	0.990&  	0.975&  	0.098&  	1&  	1&  	1&  	1\\
& 80  &0.146&  	0.997&  	0.997&  	0.997&  	0.994&  	0.102&  	1&  	1&  	1&  	1\\
& 100 &0.153&  	1    &  	0.999&  	0.999&  	0.998&  	0.098&  	1&  	1&  	1&  	1\\
\hline
  \multirow{4}{*}{Power II}
& 40  &0.635&  	0.998&  0.998& 	0.999& 	0.973& 	0.445&  	1&  	1&  	1&  	1\\
& 60  &0.850&  	1&  	1&  	1&  	0.998& 	0.703&  	1&  	1&  	1&  	1\\
& 80  &0.961&  	1&  	1&  	1&  	1&  	0.867&  	1&  	1&  	1&  	1\\
& 100 &0.992&  	1&  	1&  	1&  	1&  	0.939&  	1&  	1&  	1&  	1\\
\hline
\end{tabular}
\caption{Size and power of the test statistics evaluated under two different scenarios at 5\% significance level. Power I refers to $\mu_i \sim \frac{1}{\sqrt{p}} U(2,3)$ and Power II refers to $\mu_i \sim \frac{1}{\sqrt[4]{p}} U(2,3)$.}
\label{tab:simulation1}
\end{table}

As shown in Table \ref{tab:simulation1}, all tests for independent observations fail in controlling a given size  of $5\%$ in that these four tests
are too liberal except the independent case in Model I. On the other hand, the proposed test controls Type I error and achieves reasonable power for all four models. These results demonstrate the importance of recognizing dependence among observations in testing procedures. QQ-plots of the p-values for the four models are presented in \nameref{sec:suppl}

\subsection{Two sample case}
\label{tscase}
For the two sample case, the data was generated using the factor model \eqref{eqn:model} separately for the two groups,
\begin{eqnarray*}
 \twofam{X}{a}(t) = \twofam{\mu}{a} + \mathop{\sum}_{h=0}^{M} \twofam{A}{a}_h \twofam{\varepsilon}{a}_{t-h}, \, \, \, 1 \leq t \leq n_a, \, \, \,a = 1, 2.
\end{eqnarray*}
Autocovariance structures for the two groups can be set unequal by specifying different mixing matrices. The mixing matrices and the error covariance matrices are generated similar to the one sample case,
\begin{align}
 \twofam{A}{k}_h(i,j) = \left\{
     \begin{array}{cc}
       \frac{\twofam{\phi_1}{k}}{h |i-j|^2} & \mbox{  for  } |i-j| \leq p \twofam{w}{k}, \\
       0 & \mbox{  for  } |i-j| > p \twofam{w}{k},
     \end{array}
   \right.
   \nonumber
\end{align}
\begin{align}
   \twofam{\Sigma}{k}(i,j) = \left\{
     \begin{array}{cc}
       \frac{\sqrt{\twofam{\sigma_i}{k} \twofam{\sigma_j}{k}} \twofam{\phi_2}{k}}{|i-j|^2} & \mbox{  for  } |i-j| \leq p \twofam{w}{k}, \\
       0 & \mbox{  for  } |i-j| > p \twofam{w}{k}.
     \end{array}
   \right.
   \nonumber
\end{align}
Values of the parameters $(\twofam{w}{1} \twofam{w}{2} \twofam{\phi_1}{1}, \twofam{\phi_1}{2}, \twofam{\phi_2}{1}, \twofam{\phi_2}{2})$ are fixed as $(0.9, 0.5,$ $ 0.2, 0.4, 0.3, 0.5)$ and the two groups are taken to be balanced, with $n_1 = n_2 = n = \frac{p}{4}$. The model is tested for three levels of dependency $M = 1, 2, 3$ and three sample sizes $n = 40, 60, 80$. To study the power of the test statistic, the difference of means is generated uniformly at two rates of decay. In Model I, the difference of means is generated as $\twofam{\mu}{1} - \twofam{\mu}{2} \sim \frac{1}{\sqrt{p}} U \left(1,2 \right)$ and in Model II the difference is generated as $\twofam{\mu}{1} - \twofam{\mu}{2} \sim \frac{1}{p^{1/4}} U \left(1,2 \right)$. Type I error and power under Models I and II calculated from 10,000 randomly generated samples for all the scenarios is tabulated in Table \ref{tab:twosampleresults}. Since the proposed test statistic is shown to outperform all other test statistics for the one-sample case, results are demonstrated only for $T_{new}$.

\begin{table}[htpb]
 \centering
 \begin{tabular}{|c|c|ccc|}
 \hline
  &     & \multicolumn{3}{c|}{Sample Size} \\
  & $M$ & $40$ & $60$ & $80$ \\
  \hline
  \multirow{3}{*}{Size} 
&  1  &  0.0866  &  0.0891  &  0.0773 \\
&  2  &  0.0767  &  0.0758  &  0.0696 \\
&  3  &  0.0627  &  0.0651  &  0.0607 \\
\hline
\multirow{3}{*}{Power - Model I} 
&  1  &  0.1128  &  0.1190  &  0.1239\\
&  2  &  0.0860  &  0.0915  &  0.0882\\
&  3  &  0.0642  &  0.0734  &  0.0735\\
\hline
\multirow{3}{*}{Power - Model II} 
&  1  &  0.3003  &  0.4222  &  0.6070 \\
&  2  &  0.1712  &  0.2299  &  0.3525 \\
&  3  &  0.1103  &  0.1594  &  0.2369 \\
  \hline
 \end{tabular}
 \caption{Type I error of the proposed test statistic for the two sample case at $5\%$ significance level.}
\label{tab:twosampleresults}
\end{table}

From Table \ref{tab:twosampleresults}, it can be seen that as the sample size increases, the type I error gets closer to the significance level. From Models I and II, the power is also seen to be increasing with sample size. As expected, the power in Model II is greater than that in Model I as the norm of difference of the means is a constant in the first model whereas it increases at the rate $\sqrt{p}$ in the second model. QQ-plots of the p-values are presented in \nameref{sec:suppl}.

\subsection{Effect of Misspecification}

\begin{figure}
\centering
 \includegraphics[scale=0.35]{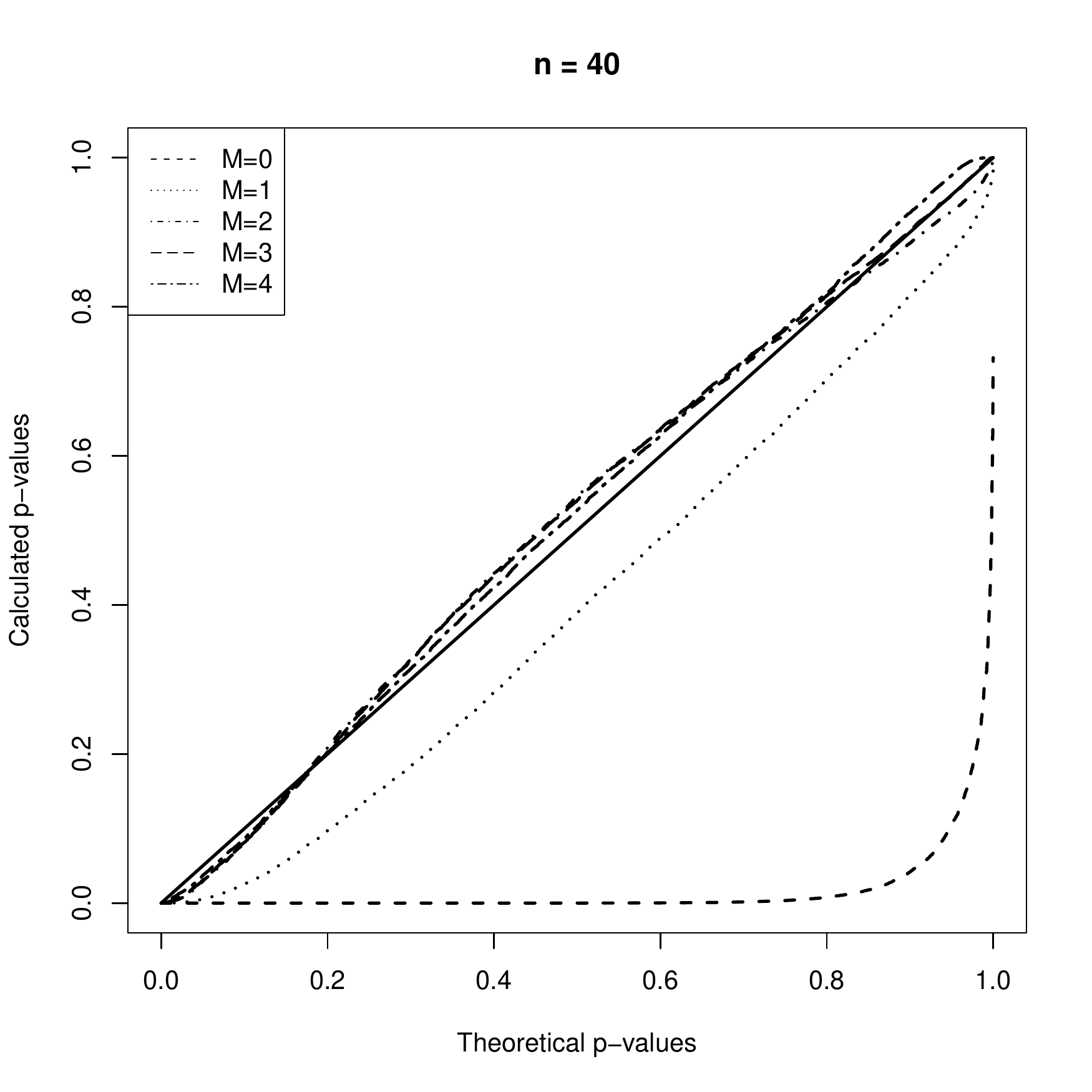} 
 \includegraphics[scale=0.35]{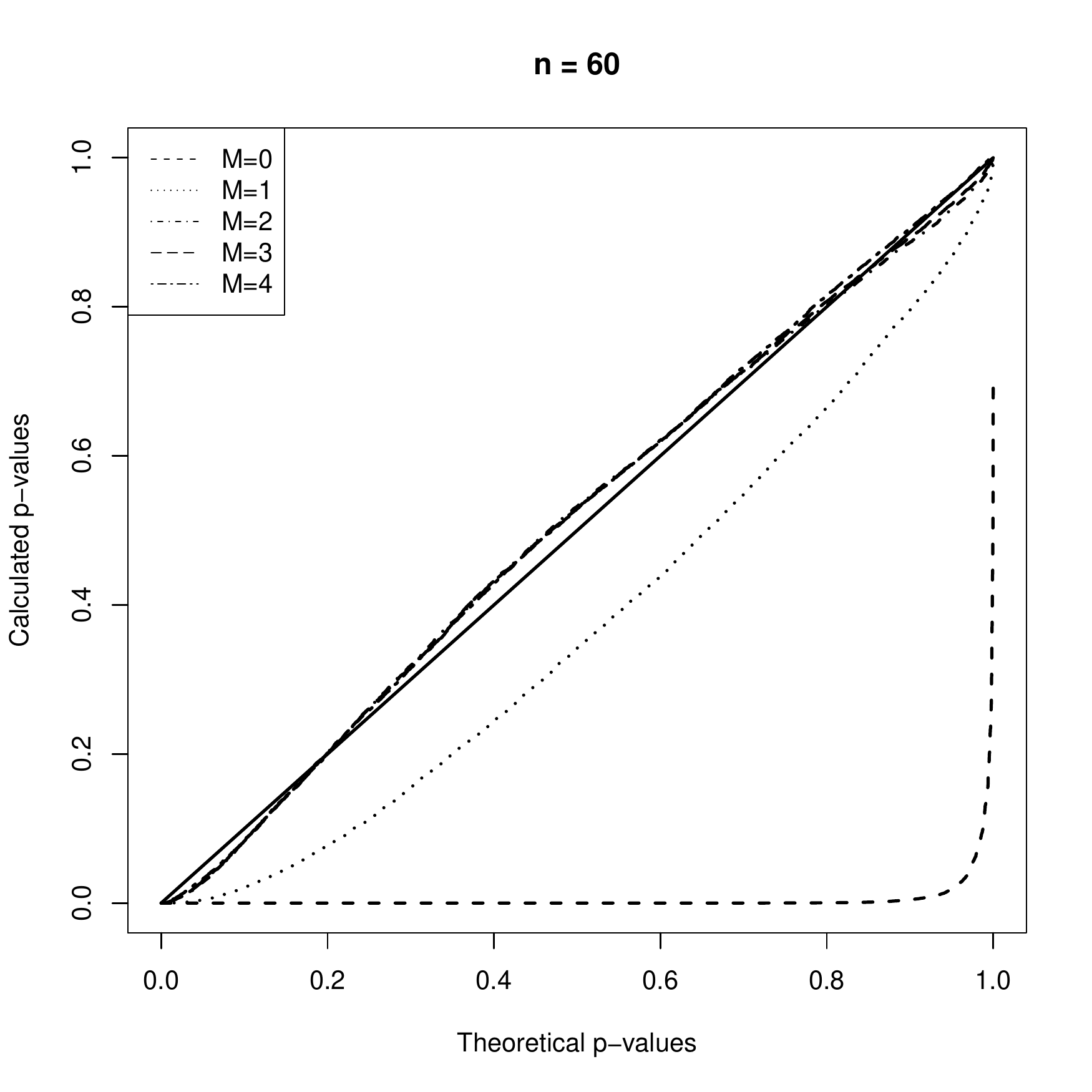} \\
 \includegraphics[scale=0.35]{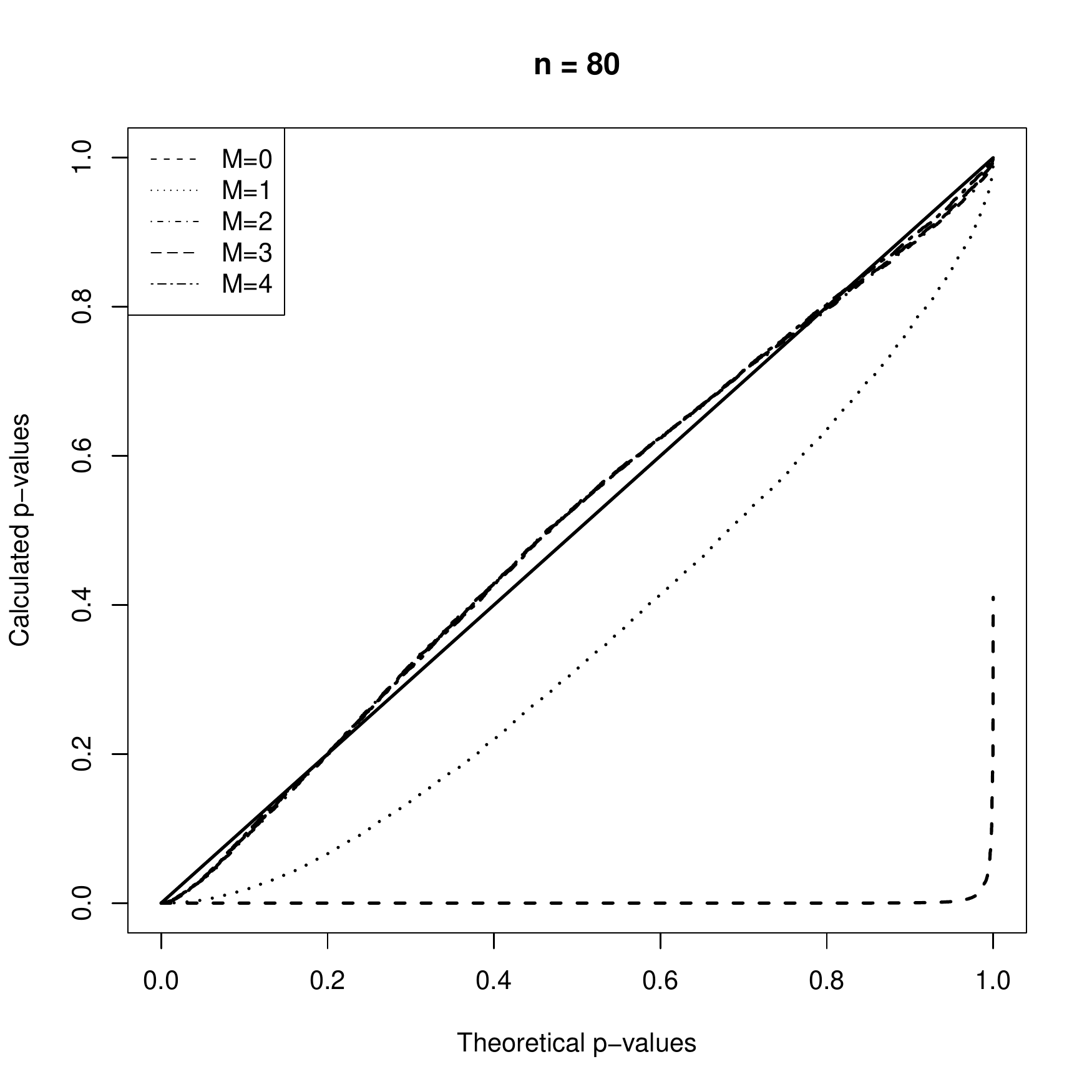} 
 \includegraphics[scale=0.35]{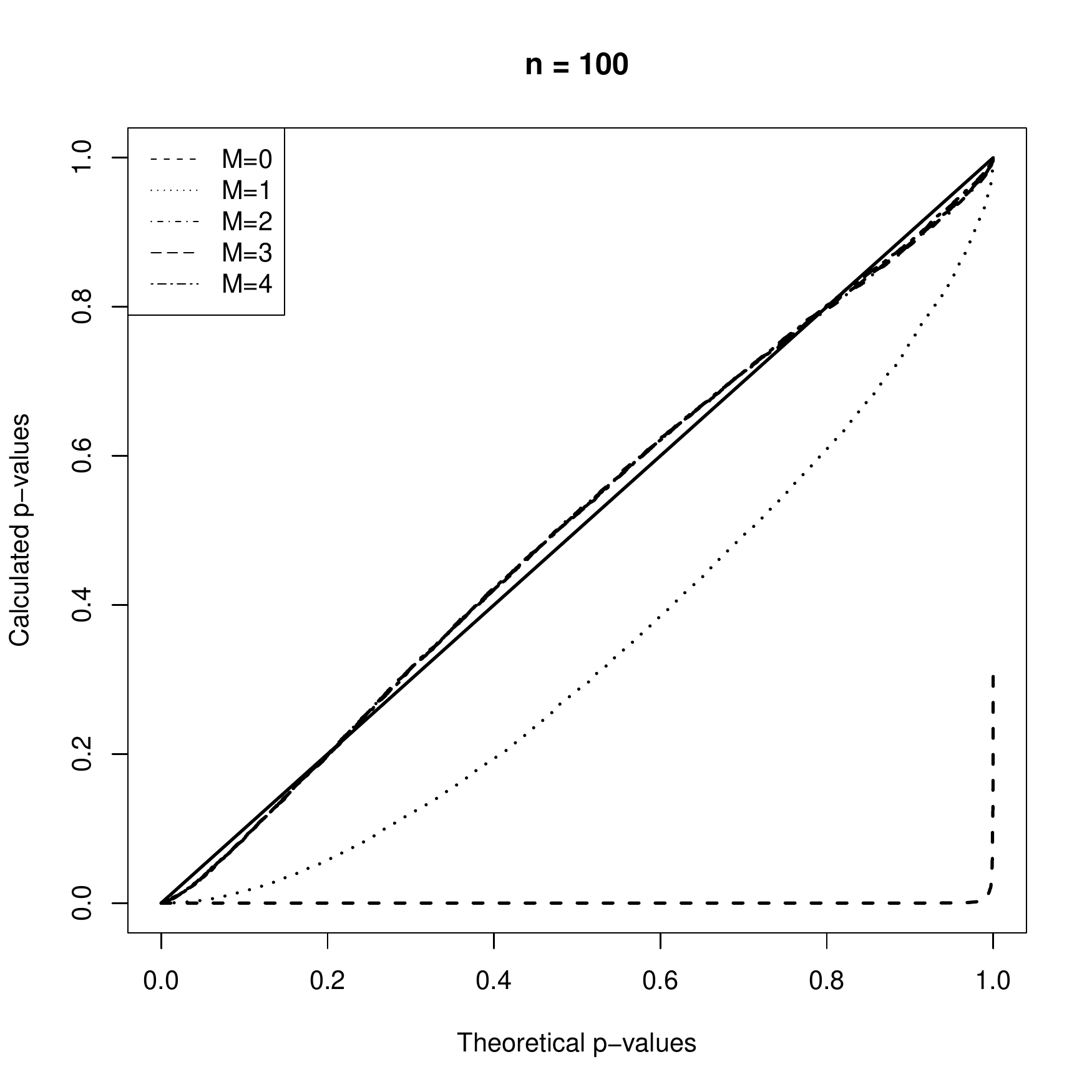}
 \caption{QQ-plot of the p-values for different levels of dependency specified and sample sizes equal to 40, 60 80 and 100 respectively.}
 \label{fig:specifym}
\end{figure}

For the simulation studies in the above sections, the value of $M$ used to simulate the data is used for calculating the test statistic. However in practical applications the value of $M$ is most likely unknown. To study the effect of misspecification of $M$ when calculating the test statistic, we tested data generated using Model III of the one sample case by specifying different values for $M$. In Model III, data was generated with the true value of $M$ equal to two and as shown in Table \ref{tab:simulation1}, $T_{new}$ is observed to outperform all the existing tests when the true value of $M$ is known. The simulation study was repeated varying the value of $M$ between 0 and 4. Type I error and power calculated using 10,000 randomly generated data sets are tabulated in Table \ref{tab:simulation2}. As seen from the table, $T_{new}$ becomes more conservative as $M$ increases, with the Type I error getting closer to the nominal size of $5 \%$. The increasing conservativeness with $M$ is also seen in the decreased power. This behavior is also persistent in the QQ-plots of the p-values, presented in Figure \ref{fig:specifym}. 

\begin{table}[htpb]
 \centering
 \begin{tabular}{cc|ccccc}
& & \multicolumn{5}{c}{Specified value of $M$} \\
&Sample size $n$ & $0$  & $1$ & $2$ & $3$ & $4$ \\
\hline
\multirow{4}{*}{Size}
 & $40$  & 0.9110  &  0.1410  &  0.0717  &  0.0685  &  0.0631 \\
 & $60$  & 0.9741  &  0.1579  &  0.0709  &  0.0715  &  0.0687 \\
 & $80$  & 0.9947  &  0.1725  &  0.0683  &  0.0669  &  0.0656 \\
 & $100$ & 0.9987  &  0.1847  &  0.0645  &  0.0647  &  0.0636 \\
 \hline
\multirow{4}{*}{Power I} 
 & $40$  & 0.9371  &  0.2178  &  0.1245  &  0.1202  &  0.1064 \\
 & $60$  & 0.9853  &  0.2598  &  0.1350  &  0.1321  &  0.1260 \\
 & $80$  & 0.9970  &  0.2965  &  0.1458  &  0.1426  &  0.1366 \\
 & $100$ & 0.9993  &  0.3232  &  0.1532  &  0.1510  &  0.1482 \\
 \hline
\multirow{4}{*}{Power II} 
 & $40$  & 0.9976  &  0.7613  &  0.6347  &  0.5995  &  0.5471 \\
 & $60$  & 1.0000  &  0.9291  &  0.8504  &  0.8382  &  0.8183 \\
 & $80$  & 1.0000  &  0.9871  &  0.9609  &  0.9576  &  0.9530 \\
 & $100$ & 1.0000  &  0.9978  &  0.9916  &  0.9904  &  0.9892 \\
\hline
\end{tabular}
\caption{Size and power of the proposed test statistic evaluated at 5\% significance level for different values of $M$ specified when true value is $M=2$. The scenarios Power I and Power II are as defined in Section \ref{sec:oscase}.}
\label{tab:simulation2}
\end{table}

\section{Conclusion}
\label{sec:conclusions}
We proposed a new test for the mean vector in high dimensional dependent data. The proposed test statistic was shown to be asymptotically normal under  conditions on the dimension, the  level of dependence and the structure of the autocovariance matrices. Through simulation studies, we illustrated that using existing tests that ignore dependence between observations leads to grossly inflated Type I error. The proposed test, being designed to incorporate dependence among the observations, reasonably controlled Type I error under a variety of dependence scenarios. When the samples were generated from an independent model($M=0$), the proposed test still performed well compared to the tests that are designed for the independent case. When the observations are assumed to be independent, $T_{new}$ does not reduce to any of the existing test statistics, since the numerator of the proposed test and the variance estimator were fashioned after $T_{BS}$ and $T_{CQ}$, respectively. A consequence of this observation is that for independent samples, $T_{BS}$ can be extended to heterogeneous covariance case by replacing the variance estimator with that from $T_{CQ}$. When the exact value of $M$ is unknown, over-specifying $M$ makes $T_{new}$ conservative. More adaptive choice of $M$ is a topic of future investigation. 

\appendix

\section{Appendix}
\label{sec:appendix}

Before constructing proofs of the theorems, we establish a few properties of the autocovariance structure which are derived from the assumptions made on the model. The following lemma is a direct extension of the model assumption stated in \eqref{eqn:condition1}.
\begin{lemma}
 \label{lem:gammaproplemma}
 For any $-M \leq a,b \leq M$ and \eqref{eqn:dim_sample},
 \begin{eqnarray*}
  \otrab{a}{b} = o \left(n^{1/2} tr \left(\Omega_n^2 \right) \right).
 \end{eqnarray*}
\end{lemma}

\noindent \textbf{Proof}: From matrix algebra, a standard result is that for any real, symmetric $p \times p$ matrix A, $tr(A) \leq \sqrt{p \, tr \left(A^2\right)}$.  Considering $b=-a$, the matrix $\Gamma(a) \Gamma(b) = \Gamma(a) \Gamma(a)'$ is symmetric and non-negative definite. Therefore,
\begin{eqnarray}
 tr \left( \Gamma(a) \Gamma(a)' \right) &\leq& \sqrt{p} \, \sqrt{tr \left( \Gamma(a) \Gamma(a)' \right)^2} \nonumber \\
  &=& \sqrt{p} \, \sqrt{tr \left( \Gamma(a) \Gamma(-a)\Gamma(a) \Gamma(-a) \right)} \nonumber \\
  &=& \sqrt{p} \, o\left(tr \left(\Omega_n^2 \right) \right), \label{eqn:prop1}
\end{eqnarray}
where the last equality is from \eqref{eqn:condition1}. For any $a$ and $b$, consider the matrix $( \Gamma(a) + \Gamma(-b) )( \Gamma(a) + \Gamma(-b) )'$. Since it is symmetric, the trace can be expressed as
 $tr [ ( \Gamma(a) + \Gamma(-b) ) ( \Gamma(a) + \Gamma(-b) )] =tr[ \Gamma(a) \Gamma(-a)] + tr [ \Gamma(b) \Gamma(-b) ] + 2 tr [ \Gamma(a) \Gamma(-b) ]$ and from \eqref{eqn:prop1}, we have
 \begin{align*}
  tr \left[ \left( \Gamma(a) + \Gamma(-b) \right) \left( \Gamma(a) + \Gamma(-b) \right)\right] &\leq  \sqrt{p} \,\sqrt{ tr \left[ \left( \Gamma(a) + \Gamma(-b) \right) \left( \Gamma(a) + \Gamma(-b) \right)\right] ^2} \\
  & =  \sqrt{p} \, \sqrt{\dmsum_{(a_1,a_2,a_3,a_4) \in D} tr \left[ \Gamma(a_1) \Gamma(a_2) \Gamma(a_3) \Gamma(a_4) \right]},
 \end{align*}
where $D = \{a,-b\} \times \{-a,b\} \times \{a,-b\} \times \{-a,b\}$ whose cardinality is finite. Therefore, it can be established that $tr \left[ \left( \Gamma(a) + \Gamma(-b) \right) \left( \Gamma(a) + \Gamma(-b) \right)\right] = \sqrt{p} \,  o\left(tr \left(\Omega_n^2 \right) \right)$, which gives
\begin{eqnarray*}
 tr \left[ \Gamma(a) \Gamma(-a) \right] + tr \left[ \Gamma(b) \Gamma(-b) \right] + 2 tr \left[ \Gamma(a) \Gamma(-b) \right] &=& \sqrt{p} \,  o\left(tr \left(\Omega_n^2 \right) \right)
\end{eqnarray*}
and by combining the above equation with equation \eqref{eqn:prop1}, we can conclude that $tr (\Gamma(a) \Gamma(b)) $ $=  \sqrt{p} o ( tr ( \Omega_n^2 ) ) = o (n^{1/2} tr (\Omega_n^2 ) )$ holds for all $a$ and $b$, where the last equality is from \eqref{eqn:dim_sample}. \qed

We now establish some matrix results which are direct extensions of the local alternative assumption for the population mean. The following property holds for any general set of matrices satisfying the conditions mentioned below while the lemma addresses the rate of increase of an inner product of $\mu$. 

\begin{property}
 \label{prop:matprop1}
 Let A and B be $n \times n$ real valued matrices such that B is symmetric and positive-definite. For all vectors $x \in \mathbb{R}^n$,
 \begin{eqnarray*}
  x' A' B A x \leq x' \left(A' A \right)^{\frac{1}{2}} x \sqrt{tr \left( A A' B^2 \right)}.
  \label{eqn:matcond1}
 \end{eqnarray*}
\end{property}

\noindent \textbf{Proof}: Let $A = U D V'$ be the singular value decomposition of the matrix $A$ where $U$ and $V$ are orthogonal matrices. This gives $A'A = V D^2 V'$, which implies $ \left(A'A \right)^{\frac{1}{2}} = V D V'$ and $A A' = U D^2 U'$. Then we have
\begin{eqnarray*}
 x' A' B A x  &=& x' V D U' B U D V' x = \left(x' V D^{\frac{1}{2}} \right) \left( D^{\frac{1}{2}} U' B U D^{\frac{1}{2}} \right) \left( D^{\frac{1}{2}} V' x \right) \\
 & \leq & \left(x' V D^{\frac{1}{2}} D^{\frac{1}{2}} V' x \right) \sqrt{tr \left( D^{\frac{1}{2}} U' B U D^{\frac{1}{2}}\right)^2 } \\
 & \leq & \left(x' V D V' x \right) \sqrt{tr \left(  U D^2 U' B^2\right)}  \\
 & = & \left(x' \left(A' A \right)^{\frac{1}{2}} x \right) \sqrt{tr \left(  A A'  B^2\right)}.
\end{eqnarray*} 
The above result uses the following preliminary matrix algebra results.
\begin{enumerate}[(i)]
 \item $|x' A x| \leq x' x \sqrt{tr \left(A' A \right)}$ for all $x \in \mathcal{R}^n$ and $A \in \mathcal{R}^{n \times n}$.
 \item If $A$ and $B$ are real symmetric matrices, then $tr \left( A B \right)^2 \leq tr \left(A^2 B^2 \right)$.
\end{enumerate}

\begin{lemma}
 \label{lem:meancond1}
 For $\mu$ which satisfy the local alternative condition mentioned in \eqref{eqn:alternativemean},
 \begin{eqnarray*}
  \mu' \Gamma(-a) \Gamma(0) \Gamma(a) \mu = o \left( n^{-1} tr^2 \left( \Omega_n^2 \right) \right), \, \, \forall -M \leq a \leq M.
 \end{eqnarray*}
\end{lemma}

\noindent \textbf{Proof}: From Property \eqref{prop:matprop1}, we have
\begin{eqnarray*}
 \frac{\mu' \Gamma(-a) \Gamma(0) \Gamma(a) \mu}{tr^2 \left( \Omega_n^2 \right)} & \leq & \frac{ \left(\mu' \left(\Gamma(a) \Gamma(-a) \right)^{\frac{1}{2}} \mu \right) \sqrt{tr \left(\Gamma(a) \Gamma(-a) \Gamma(0)^2 \right)}}{tr^2 \left( \Omega_n^2 \right)} \\
 & = & o \left(n^{-1} \right).
\end{eqnarray*}
where the rate is due to \eqref{eqn:alternativemean} and \eqref{eqn:condition1}.

\begin{lemma}
 \label{lem:alternative2}
 For any $\mu$ which satisfy the local alternative condition in \eqref{eqn:alternativemean} and lag $-M \leq a \leq M$,
 \begin{eqnarray}
  \mu' \Gamma(a) \mu = o \left( n^{-1} tr \left(\Omega_n^2 \right) \right).
  \label{eqn:alternativemean2}
 \end{eqnarray}
\end{lemma}

\noindent \textbf{Proof}:
The autocovariance matrix at lag zero is symmetric and positive definite. Therefore the result follows immediately for $a=0$. For non-zero $a$, Let $\Gamma(a) = U D V'$ be the singular value decomposition, where $U$ and $V$ are orthogonal matrices and $D$ is diagonal with positive entries. Therefore, we have
\begin{eqnarray*}
 \left| \mu' \Gamma(a) \mu \right| &=& \left|\mu' U D V' \mu \right| =\left| \left[ \mu' U D^{\frac{1}{2}} \right]\left[  D^{\frac{1}{2}} V' \mu \right]\right|  \\
  & \leq & \sqrt{\left( \mu' U D U' \mu \right) \left( \mu' V D V' \mu \right)} \\
  & \leq & \sqrt{ \left( \mu' \left( A A' \right)^{\frac{1}{2}} \mu \right)  \left( \mu' \left( A' A \right)^{\frac{1}{2}} \mu \right)}
  = o \left( n^{-1} tr \left(\Omega_n^2 \right) \right),
\end{eqnarray*}
where $A = U D^{1/2}V'$. Note however that \eqref{eqn:alternativemean2} does not imply the local alternative condition.

In the following subsections, we shall establish asymptotic normality of the test statistic and prove ratio consistency of the variance estimator.

\subsection{Proof of Theorem  \ref{thm:clt}}

Let $X_i, i = 1, 2, \ldots, n$ be a realization of a $M$-dependent stationary Gaussian process. We shall establish asymptotic normality of the test statistic under the null hypothesis. The proof follows by argument of triangular arrays which has been extensively in literature \cite{Diananda,Hoeffding}. Central Limit Theorem for $M$-dependent univariate variables has been discussed in \cite{Berk,RomanoWolf}. Since the test statistic involves inner product of data vectors, we propose an extension of the univariate triangular array to a two dimensional array. To establish asymptotic normality, express the test statistic as follows,
\begin{eqnarray}
 \mtxt{T} = \frac{\mtxt{M}}{\sqrt{Var \left(\mtxt{M}\right)}} &=& \frac{\overline{X}' \overline{X} - \widehat{tr \left(\Omega_n\right)}}{\sqrt{Var \left(\mtxt{M}\right)}}, \nonumber \\
 &=& \frac{\overline{X}' \overline{X} - tr \left(\Omega_n\right)}{\sqrt{Var \left(\mtxt{M}\right)}} - \frac{\widehat{tr \left(\Omega_n\right)} - tr \left(\Omega_n\right)}{\sqrt{Var \left(\mtxt{M}\right)}}, \nonumber \\
 &=& \Delta_1 - \Delta_2,
 \label{eqn:delta12}
\end{eqnarray}
where the two terms are defined as $\Delta_1 = \frac{\overline{X}' \overline{X} - tr \left(\Omega_n\right)}{\sqrt{Var \left(\mtxt{M}\right)}}$ and $\Delta_2 =\frac{\widehat{tr \left(\Omega_n\right)} - tr \left(\Omega_n\right)}{\sqrt{Var \left(\mtxt{M}\right)}} $. We shall establish asymptotic normality of the first term using two dimensional triangular arrays. The proof of asymptotic normality of $\mtxt{T}$ is completed by showing that the second term converges to zero in probability.

\begin{figure}[htp]
\centering
\includegraphics[width=1.0\textwidth]{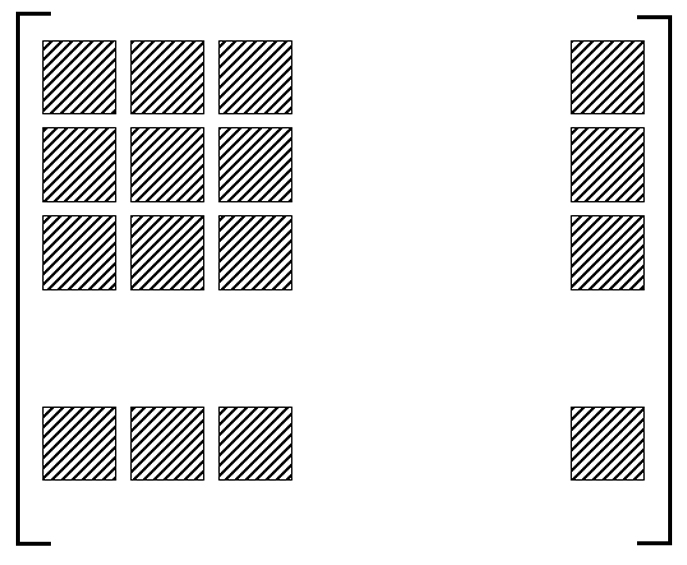}
\caption{Pictorial representation of the two dimensional triangular array.}
\label{fig:array}
\end{figure}

Let $\mathcal{A}$ be a $n \times n$ matrix with elements given by $\mathcal{A}_{ij} = \frac{1}{n^2} ( X_i'X_j - tr \left(\Gamma \left(|i-j|\right)\right) )$ for $1 \leq i,j \leq n$. The sum of all elements of $\mathcal{A}$ is equal to $\overline{X}' \overline{X} - tr \left(\Omega_n\right)$, which is the numerator of $\Delta_1$. This motivates using the matrix $\mathcal{A}$ as our two dimensional triangular array. For any n, choose $0 < \alpha < 1$, $C > 0$ and $w_n = C n^{\alpha} > M$, so that
\begin{eqnarray}
 n = w_n k_n + r_n.
 \label{eqn:wnkn}
\end{eqnarray}
The quantity $w_n$ represents the size of each block to be constructed, $k_n$ represents the number of blocks constructed in each direction and $r_n$ represents the dimension of the remainder terms. Define the random variables
\begin{eqnarray}
 B_{ij} &=& \msum_{k = (i-1)w_n + 1}^{i w_n - M} \msum_{l = (j-1)w_n + 1}^{j w_n - M} \mathcal{A}_{kl} \nonumber \\
 D_{ij} &=& \left(\msum_{k = (i-1)w_n + 1}^{i w_n} \msum_{l = (j-1)w_n + 1}^{j w_n} \mathcal{A}_{kl} \right) - B{ij} \nonumber \\
 F &=& \dmsum_{\left(i,j\right) \in \{1,2,\ldots,n\}^2 \smallsetminus \{1,2,\ldots,w_n k_n\}^2} \mathcal{A}_{ij}  \label{eqn:triangarray}
\end{eqnarray}
for $1 \leq i,j \leq k_n$. Since the matrix $\mathcal{A}$ is symmetric, $B_{ij} = B_{ji}$. The random variables $B_{ij}$'s are constructed from blocks in the matrix $\mathcal{A}$ which are at least $M$ indices apart. Pictorial representation of the matrix $\mathcal{A}$ describing the construction of the blocks is provided in Figure \ref{fig:array}. The shaded blocks represent the variables $B_{ij}$'s and the unshaded portion corresponds to $D_{ij}$'s and $F$. The shaded blocks are of dimension $ \left(w_n - M \right) \times \left(w_n - M \right)$ and the blocks are exactly $M$ indices apart. For an appropriate choice of $\alpha$, as the sample size increases, the shaded region becomes dominant with respect to the unshaded region. Hence the contribution of the terms coming from the unshaded region will be negligible in comparison to the terms coming from the shaded region. 

Since the process is Gaussian with M-dependent stationarity, $ \{B_{ij} , 1 \leq i,j \leq k_n \}$ are independent. By construction, we have $E \left(B_{ij} \right) = 0$ for all $1 \leq i, j \leq k_n$. However, they are not identical variables because their variance depends on the index,
\begin{align}
Var \left(B_{ij} \right) = \frac{(w_n - M)^2\left(1 + \mathcal{J}(i-j,0) \right)}{n^4}   tr \left( \msum_{h=-M}^M \left(1 - \frac{|h|}{w_n-M} \right) \Gamma(h) \right)^2,
\end{align}
where $\mathcal{J}$ is the indicator function, with $\mathcal{J}(a,b) = 1$ if $a = b$ and $\mathcal{J}(a,b) = 0$ if $a \neq b$. From the above constructed random variables, construct a sequence of random variables $\left\{Y_i\right\}, 1 \leq i \leq \frac{k_n(k_n+1)}{2}$ which is given by
\begin{eqnarray}
 Y_{(i-1)k_n - \frac{(i-1)(i-2)}{2} + j} = \left(2 - \mathcal{J}\left(i-j,0\right) \right) B_{ij}.
 \label{eqn:cltrvs}
\end{eqnarray}
This sequence $\left\{Y_i\right\}$ is the vector form of the matrix $\left(B_{ij}\right){1 \leq i,j \leq k_n}$ and the coefficients are to adjust for multiplicity.
In the above construction, the off-diagonal blocks are multiplied by two because of symmetry of $\mathcal{A}$. 

To establish asymptotic normality of $\Delta_1$, express the quantity as
\begin{eqnarray}
 \Delta_1 &=& \frac{\overline{X}' \overline{X} - tr \left(\Omega_n \right)}{\sqrt{\varmn}} = \frac{\dmsum_{i,j=1}^n \mathcal{A}_{ij}}{\sqrt{\varmn}} \nonumber\\
 &=& \frac{ \dmsum_{i,j=1}^{k_n} B_{ij}}{\sqrt{\varmn}} + \frac{  \dmsum_{i,j=1}^{k_n} D_{ij} + F}{\sqrt{\varmn}} \nonumber \\
 &=& \Delta_{11} + \Delta_{12}. \label{eqn:delta1}
\end{eqnarray}
By construction, $\dmsum_{i,j=1}^{k_n} B_{ij} = \sum_{i=1}^{k_n(k_n-1)/2} Y_i$ where $\{Y_i\}$ is a sequence of independent variables. The sum of variance of $Y_i$'s is 
\begin{eqnarray*}
 \msum_{i=1}^{\frac{k_n(k_n+1)}{2}} Var\left(Y_i \right) &=& k_n Var \left(B_{11} \right) + \frac{k_n \left(k_n-1\right)}{2} Var \left(2 B_{12} \right) \\
 &=& k_n \frac{2}{n^2}  \left( \frac{w_n - M}{n} \right)^2 tr \left( \msum_{h=-M}^M \left(1 - \frac{|h|}{w_n-M}\right) \Gamma(h)  \right)^2 \\
 &+&  4 \frac{k_n \left(k_n-1\right)}{2} \frac{1}{n^2}  \left( \frac{w_n - M}{n} \right)^2 tr \left( \msum_{h=-M}^M \left(1 - \frac{|h|}{w_n-M}  \right) \Gamma(h) \right)^2 \\
 &=& k_n^2 \frac{2}{n^2}  \left( \frac{w_n - M}{n} \right)^2 tr \left( \msum_{h=-M}^M \left(1 - \frac{|h|}{w_n-M}\right) \Gamma(h)  \right)^2. \\
\end{eqnarray*}
By construction of $w_n$, $k_n$ and $r_n$, we have $\frac{\msum_{i=1}^{\frac{k_n(k_n+1)}{2}} Var\left(Y_i \right)}{\varmn} \rightarrow 1$. Therefore asymptotic normality of $\frac{ \dmsum_{i,j=1}^{k_n} B_{ij}}{\sqrt{\varmn}}$ holds if the Lyapunov condition is satisfied. To verify the Lyapunov condition, the fourth order central moment of $B_{ij}$'s are
\begin{eqnarray*}
 E \left(B_{ii}^4\right) &=& \left( \frac{w_n-M}{n} \right)^4 \left[ 12 tr^2 \left(\Omega_{w_n}^2 \right) + 48 tr \left(\Omega_{w_n}^4 \right) \right], \\
 E \left(B_{ij}^4\right) &=& 6 \left( \frac{w_n-M}{n} \right)^4 \left[  tr^2 \left(\Omega_{w_n}^2 \right) +  tr \left(\Omega_{w_n}^4 \right) \right],
\end{eqnarray*}
where $\Omega_{w_n}$ is as defined in \eqref{eqn:omega}. Therefore we have
\begin{align}
\frac{\msum_{i=1}^{\frac{k_n(k_n+1)}{2}} E \left(Y_i^4\right)}{\varmn^2}  =  \frac{12}{n^4} \left( \frac{w_n-M}{n} \right)^4 \frac{ \left(4 k_n^2 - 3 k_n \right) tr^2 \left(\Omega_{w_n}^2 \right) + 4 k_n^2 tr \left(\Omega_{w_n}^4 \right)}{\frac{4}{n^4} tr^2 \left(\Omega_n^2 \right)},
\label{eqn:lyapunov}
\end{align}
which converges to zero because of choice of $k_n$ and $w_n$ and because $\frac{tr \left(\Omega_{w_n}^4 \right)}{tr^2 \left(\Omega_n^2 \right)} = O(1)$ by Kashin-Garnaev-Gluskin inequality \cite{KGG}. 

The second term in $\Delta_1$, denoted by $\Delta_{12}$ has expected value equal to zero by construction. By symmetry of the matrix $\mathcal{A}$, we can express $\Delta_{12}$ as
\begin{eqnarray*}
 \Delta_{12} = \frac{1}{\sqrt{\varmn}} \msum_{a=1}^{k_n} \underbrace{ \left( \msum_{i=1}^n \msum_{j=a w_n - M +1}^{a w_n} \mathcal{A}_{ij} \right)}_{\equiv \Delta_{12a}} = \frac{1}{\sqrt{\varmn}} \msum_{a=1}^{k_n} \Delta_{12a} .
\end{eqnarray*}
For $a \neq b$, the terms $\Delta_{12a}$ and $\Delta_{12b}$ are independent because the width $w_n$ is greater than $M$ and increases faster with respect to $n$. By construction of the indexing sets for the summation in $\Delta_{12a}$s, they have equal variance. Therefore the variance of $\Delta_{12}$ can be expressed as
\begin{eqnarray}
 Var \left( \Delta_{12} \right) &=& Var \left[\frac{1}{\sqrt{\varmn}} \msum_{a=1}^{k_n} \left( \msum_{i=1}^n \msum_{j=a w_n - M +1}^{a w_n} \mathcal{A}_{ij} \right)\right] \nonumber \\
 &=& Var \left[\frac{1}{\sqrt{\varmn}} \msum_{a=1}^{k_n} \Delta_{12a} \right] =   \frac{1}{\varmn} k_n Var \left(\Delta_{121} \right) \nonumber \\
 &=& \frac{1}{\varmn} \msum_{i=1}^n \msum_{j=w_n-M+1}^{w_n} \mathcal{A}_{ij} \nonumber \\
 &=& \frac{k_n}{\varmn} \frac{1}{n^4} \left[ tr \left( \msum_{i=1}^n \msum_{j=1}^M \Gamma(i-j) \right)^2 + tr \left( \msum_{i=1}^n \msum_{j=w_n-M+1}^{w_n} \Gamma(i-j) \right)^2 \right] \nonumber \\
 &=& \frac{k_n} {n^4 \varmn} \left[ n M tr \left(\Omega_n \Omega_M \right) + M^2 tr \left( \msum_{a=-M}^M \Gamma(a) \right)^2 \right] \nonumber \\
 &=& \frac{k_n M}{2n} \left( \frac{tr \left(\Omega_n \Omega_M \right)}{tr \left(\Omega_n^2 \right)} + \frac{M}{n} \frac{ tr \left( \msum_{a=-M}^M \Gamma(a) \right)^2}{tr \left( \Omega_n^2 \right)} \right) \left(1 + o(1) \right). \nonumber
\end{eqnarray}
From assumptions \eqref{eqn:M} and \eqref{eqn:wnkn} on the rate of increase of $w_n$ and $M$ with respect to $n$, we have $\frac{k_nM}{n}$ and $\frac{M}{n}$ converging to zero, while the remaining terms are finitely bounded. Therefore $Var \left(\Delta_{12} \right)$ converges to zero, which implies $\Delta_{12}$ converges to zero in probability. Hence by Slutksy's theorem, we can conclude that $\Delta_1$ is asymptotically normal.

The second term of $\mtxt{T}$ in \eqref{eqn:delta12} is equal to
\begin{eqnarray*}
 \Delta_2 = \frac{\widehat{tr \left(\Omega_n\right)} - tr \left(\Omega_n \right)}{\sqrt{\varmn}},
\end{eqnarray*}
which has expected value zero because the estimator is constructed to be unbiased. The variance is given by
\begin{eqnarray*}
 Var \left(\Delta_2 \right) &=& \frac{ Var \left(\widehat{tr \left(\Omega_n\right)} \right)}{\varmn} = \frac{Var \left( \msum_{h=0}^M \beta_n(h) tr \left(\widehat{\Gamma(h)} \right) \right)}{\varmn} \\
 &=& \frac{\dmsum_{a,b=-M}^M tr \left(Q D_a Q D_b \right) tr \left(\Gamma(a) \Gamma(b) \right)}{\varmn}
\end{eqnarray*}
where the coefficients $\left\{\beta_n(h) \right\}_{0 \leq h \leq M}$ are obtained from the expression for $\widehat{tr \left(\Omega_n \right)}$ and
\begin{eqnarray*}
 Q = \frac{1}{n} \msum_{h = 0}^M \beta_n(h) \left[D_h - \frac{1}{n} \msum_{i=0}^{n-h} e_i 1' - \frac{1}{n} \msum_{j=h+1}^n e_j 1' + \frac{n-h}{n^2} \mathcal{J} \right].
\end{eqnarray*}
The matrix $D_a$ is a $n \times n$ matrix with ones along the $a$-th diagonal and all remaining elements equal to zero and $\mathcal{J}$ is the all-ones matrix. Using the expression in \eqref{eqn:thetas} to calculate $\left\{\beta_n(h) \right\}_{0 \leq h \leq M}$, it is easy to verify that
\begin{eqnarray*}
 \max_{-M \leq a,b \leq M} \left|tr \left(Q D_a Q D_b \right) \right| = o \left(n^{-\frac{11}{4}} \right).
\end{eqnarray*}

Therefore we have
\begin{eqnarray*}
 Var \left(\Delta_2 \right) \leq \max \left|tr \left(Q D_a Q D_b \right) \right| \dmsum_{a,b=-M}^M \frac{ \left| tr \left(\Gamma(a) \Gamma(b) \right) \right|}{\frac{2}{n^2} tr \left(\Omega_n^2 \right)} \left(1 + o(1) \right), \\
\end{eqnarray*}
and by Lemma \eqref{lem:gammaproplemma} the above quantity converges to zero. Therefore $\Delta_2$ converges to zero in probability and by Slutsky's theorem, we have asymptotic normality of $\mtxt{T}$.

\qed

\subsection{{Proof of Proposition } \ref{prop:simplexhis}}

From the expression in \eqref{eqn:mn}, the quantity $\mtxt{M}$ can be expressed in a quadratic form as $\mtxt{M} = Z' \left(\Pi \otimes \mathcal{I}_p \right) Z$, where $\Pi(i,j) = \pi_n(i,j)$ and $Z = \begin{bmatrix}X_1'& X_2' &\ldots& X_n'\end{bmatrix}'$ has a multivariate normal distribution with mean zero and variance $\Sigma$ given by $\Sigma = \msum_{h = -M}^M D_h \otimes \Gamma(h)$, where $D_h$ is the matrix with ones on its $h^{th}$ diagonal and zero elsewhere, i.e. $D_h(i,j) = \mathcal{J}(i-j,h)$. This quadratic form implies the variance of $\mtxt{M}$ is
\begin{eqnarray*}
 Var \left(\mtxt{M} \right) = tr \left( \left(\Pi \otimes \mathcal{I}_p \right) \Sigma \right)^2 = \msum_{a=-M}^M \msum_{b=-M}^M tr \left(\Pi D_a \Pi D_b \right) tr \left(\Gamma(a) \Gamma(b) \right).
\end{eqnarray*}
Readjusting terms in the expression of $\pi_n(t,s)$ gives
\begin{eqnarray*}
 \Pi &=& \frac{1}{n^2} \left(1 - \chi_n \right) J - \frac{1}{2n} \msum_{h=0}^M \beta_n(h) \left( D_h + D_{-h} \right) \\
   && + \frac{1}{2 n^2} \msum_{h=0}^M \beta_n(h) \left( D_h + D_{-h} \right) J + \frac{1}{2 n^2} \msum_{h=0}^M \beta_n(h) J \left( D_h + D_{-h} \right),
\end{eqnarray*}
where $J$ is the matrix of all ones, $\beta_n = \left( \beta_n(h) \right)_{0 \leq h \leq M} = \Theta_n^{-1}b$ is the vector of coefficients obtained using \eqref{eqn:thetas} and $\chi_n = \msum_{h=0}^M \left(1 - \frac{h}{n} \right) \beta_n(h)$. Since $M = O \left(n^{\frac{1}{8}} \right)$ , the number of terms in $\chi_n$ increases. As $n$ approaches infinity, $\chi_n$ becomes an infinite summation. The elements of $\beta_n$ are of order $\frac{1}{n}$, therefore $\chi_n$ converges to zero as $n$ goes to infinity. Evaluating and simplifying the terms, the given expression for $\xi_n(a,b)$ can be derived. Replacing the coefficients by the corresponding approximating quantities in the expression of variance,
\begin{eqnarray*}
 \varmn &=& 2 \msum_{a=-M}^M \msum_{b=-M}^M \xi_n(a,b) tr \left( \Gamma(a) \Gamma(b) \right) \\
  &=& \frac{2}{n^2} \msum_{a=-M}^M \msum_{b=-M}^M  \left(1 + \chi_n \right) \left( 1 - \frac{|a|}{n} \right)  \left( 1 - \frac{|b|}{n} \right) tr \left( \Gamma(a) \Gamma(b) \right) \left(1 + o(1) \right) \\
  &=& \frac{2}{n^2} \left(1 + \chi_n \right) tr \left[ \msum_{a = -M}^M \Gamma(a) \right]^2 \left(1 + o(1) \right) \\
  &=& \frac{2}{n^2} \left(1 + \chi_n \right) tr \left( \Omega_n^2 \right) \left(1 + o(1) \right).
\end{eqnarray*}
Therefore, it follows that $\frac{\varmn}{2 n^{-2} tr \left(\Omega_n^2 \right)} =  \left(1 + \chi_n \right)\left(1 + o(1) \right) \stackrel{n \rightarrow \infty}{\longrightarrow} 1$.
\qed

The following lemma \eqref{lem:c1c2} and \eqref{lem:btslemma} are instrumental in proving Theorem \eqref{thm:variance}. 
\begin{lemma}
 \label{lem:c1c2}
 For any $-M \leq a,b \leq M$, define the set 
 \begin{eqnarray}
 \mathcal{A} = \left\{\left(t,s\right) : 1 \leq t,s \leq n, |t+a-s-b| > M, |t-s|>M \right\},
 \label{eqn:indexA}
 \end{eqnarray}
 and $n_{a,b}$ be its cardinality. Define $\mathcal{C}$ as
 \begin{align*}
  \mathcal{C} = \left\{\left(q_1,q_2,q_3,q_4 \right) = f\left(t_1-t_2,t_1-s_1,t_1-s_2,t_2-s_1,t_2-s_2,s_1-s_2 \right)\right\}
 \end{align*}
where $\left(t_1,s_1 \right), \left(t_2,s_2 \right) \in \mathcal{A}$ and $f : \mathbb{Z}^6 \rightarrow \mathbb{Z}^4$. The number of elements of $\mathcal{C}$ for which $\left(q_1,q_2,q_3,q_4\right) = \left(w_1,w_2,w_3,w_4\right)$ where $-M \leq w_i \leq M$ is of the same order of $n_{a,b}$. If at least one of $q_i$'s is fixed, then the number of elements will be $O \left(n^3 \right)$. 
\end{lemma}

\noindent \textbf{Proof}: We present the result for $f (t_1-t_2,t_1-s_1,t_1-s_2,t_2-s_1,t_2-s_2,s_1-s_2 ) = (t_1-t_2,s_1-s_2,t_1-t_2,s_1-s_2)$. Similar bound can be claimed for other functions. When $\left(q_1,q_2,q_3,q_4\right) = \left(t_1-t_2,s_1-s_2,t_1-t_2,s_1-s_2\right)$,
\begin{eqnarray*}
 && \dmsum_{(t_1,s_1),(t_2,s_2) \in \mathcal{A}} \mathcal{J}\left(q_1 = w_1  \right) \mathcal{J}\left(q_2 = w_2 \right) \mathcal{J}\left(q_3 = w_3 \right) \mathcal{J}\left(q_4 = w_4 \right) \\
 && = \dmsum_{(t_1,s_1),(t_2,s_2) \in \mathcal{A}} \mathcal{J}\left(t_1-t_2 = w_1 = w_3 \right) \mathcal{J}\left(s_1-s_2 = w_2 = w_4\right)  \\
 && = \dmsum_{\substack{(t_1,s_1) \in \mathcal{A}}} \mathcal{J} \left(\left(t_1-w_1, s_1-w_2\right) \in \mathcal{A} \right) \mathcal{J} \left(w_1=w_3\right)  \mathcal{J} \left(w_2=w_4 \right)\\
 && = O \left(n_{a,b} \right) = O \left(n^2 \right).
\end{eqnarray*}
The above bound holds only when all the variates in the function $f$ depend on differences $t_i - s_j$. Restricting all $q_i$'s as functions of $\left(t_1,t_2,s_1,s_2\right)$ imposes restrictions on all the four quantities. When one of the $q_i$'s is fixed, one of $t_1,s_1,t_2,s_2$ will be unrestricted. Therefore the number of pairs $\left(t_1,s_1\right),\left(t_2,s_2\right)$ which contribute to the summation will be of the order $n_{a,b} \, n = O \left(n^3 \right)$.
\qed 

For any $-M \leq a,b \leq M$ and $(t,s) \in \mathcal{A}(a,b)$, recall the definition of the set $\mathcal{B}(t,s) = \{ 1 \leq i \leq n : |i-t|>M, |i-s| > M, |i-t-a|>M, |i-s-b|>M \}.$ This set consists of all the indices such that the observations $\left\{X_w\right\}_{w \in \mathcal{B}(t,s)}$ are independent of $\left(X_t,X_s,X_{t+a},X_{s+b}\right)$. The following lemma establishes the rate of variance of the average of observations in the set. 
\begin{lemma}
 \label{lem:btslemma}
  The variance of the inner product ${\overline{X}_n^*}' \overline{X}_n^*$, where $\overline{X}_n^*$ is the average of $\{X_w , w \in \mathcal{B}(t,s) \}$  is of the same order as $\varmn$, that is $Var \left( {\overline{X}_n^*}' \overline{X}_n^* \right) = \frac{2}{m_{\mathcal{B}}^2} tr \left( \Omega_{n \mathcal{B}}^{*2}\right)$ where
 \begin{eqnarray*}
 \Omega_{n \mathcal{B}}^*   = \msum_{h = -M}^M a_{n,\mathcal{B}}(h) \Gamma(h),
 \end{eqnarray*}
 $m_{\mathcal{B}}$ is the cardinality of $\mathcal{B}$ with $0 <  lim inf \frac{m_{\mathcal{B}}}{n} \leq lim sup \frac{m_{\mathcal{B}}}{n} < 1$ and the coefficients are of the same order as those in \eqref{eqn:omega},  $a_{n,\mathcal{B}}(h) = \left(1 - \frac{|h|}{n} \right) \left(1 + o(1) \right)$ for $-M \leq h \leq M$.
\end{lemma}

\noindent \textbf{Proof:} We prove the lemma for the case $a = b = 0$, $t = 1$ and $s = n$ and the proof can be extended for other values similarly. When $\left\{ a,b,t,s\right\} = \left\{0,0,1,n\right\}$, we have 
\begin{align*}
\mathcal{B}(t,s) &= \left\{ 1 \leq i \leq n : |i-1|>M, |i-n| > M, |i-1|>M, |i-n| > M \right\} \\
&= \left\{ M+2, M+3, \ldots, n-M-1 \right\}.
\end{align*}
Therefore $ \overline{X}_n^*$ is the average of $n-2(M+1)$ observations, which gives
\begin{eqnarray*}
\Omega_{n \mathcal{B}}^* &= & Var \left(\overline{X}_{n-2(M+1)}^{'} \overline{X}_{n-2(M+1)} \right)  \\
&=& \frac{2}{\left(n-2(M+1) \right)^2} tr \left[\mathop{\sum}_{h=-M}^M \left(1 - \frac{|h|}{n-2(M+1)} \right)\Gamma(h)  \right]^2.
\end{eqnarray*}
In this illustration, $m_{\mathcal{B}} = n - 2(M+1)$ and $a_{n, \mathcal{B}}(h) = \left(1 - \frac{|h|}{n-2(M+1)} \right)$ and it is straightforward to verify the conditions mentioned in the statement of the lemma hold. \qed

\subsection{Proof of Theorem \ref{thm:variance}}
\label{sec:varianceproof}

In this section, we derive ratio consistency of the estimator of variance of $\mtxt{M}$. The results are derived for any non-zero value of $\mu$ which satisfy the local alternative condition in \eqref{eqn:alternativemean}. 
For any value of the population mean $\mu$, let $Y_t = X_t - \mu$ be the observations centered at zero. The estimator of trace of product of two autocovariance matrices is given by
\begin{eqnarray*}
 \trab{a}{b} &=& \frac{1}{n_{a,b}} \dmsum_{t,s \in \mathcal{A}}  \left(X_{t+a} - \overline{X}_n^* \right)' X_{s} \left(X_{s+b} - \overline{X}_n^* \right)' X_{t} \\
 &=& \frac{1}{n_{a,b}} \dmsum_{t,s \in \mathcal{A}} \left(Y_{t+a} - \overline{Y}_n^* \right)' \left(Y_{s} + \mu \right) \left(Y_{s+b} - \overline{Y}_n^* \right)'.
\end{eqnarray*}
The estimator proposed for $\varmn$ can then be expressed in terms of $Y_t$ and $\mu$ as
\begin{eqnarray*}
 \widehat{\varmn} &=& \dmsum_{a,b = -M}^M \xi_n(a,b) \trab{a}{b} \\
 &=& \dmsum_{a,b = -M}^M \frac{\xi_n(a,b)}{n_{a,b}} \dmsum_{t,s \in \mathcal{A}} \left(Y_{t+a} - \overline{Y}_n^* \right)' \left(Y_{s} + \mu \right) \left(Y_{s+b} - \overline{Y}_n^* \right)' \left( Y_{t} + \mu \right) \\
 &=& \dmsum_{a,b = -M}^M  \frac{\xi_n(a,b)}{n_{a,b}} \dmsum_{t,s \in \mathcal{A}} \bigg\{ \left(Y_{t+a} - \overline{Y}_n^* \right)' Y_s \left(Y_{s+b} - \overline{Y}_n^* \right)' Y_t \\
 && +  \left(Y_{t+a} - \overline{Y}_n^* \right)' \mu \left(Y_{s+b} - \overline{Y}_n^* \right)' Y_t  \\
 && + \left(Y_{t+a} - \overline{Y}_n^* \right)' Y_s \left(Y_{s+b} - \overline{Y}_n^* \right)' \mu \\
 && + \left(Y_{t+a} - \overline{Y}_n^* \right)' \mu \left(Y_{s+b} - \overline{Y}_n^* \right)' \mu \bigg\} \\
 &=& \msum_{k=1}^4 W_k.
\end{eqnarray*}

The terms $W_1, W_2, W_3, W_4$ as expressed above will be studied individually. The first term $W_1$ consists of fourth order moments while the remaining terms consist of lower order moments. To establish ratio consistency of the proposed estimate, we shall show that
\begin{enumerate}[(i)]
 \item  $W_1 = \varmn \left(1 + o_p(1) \right)$
 \item  $W_k = o_p \left(\varmn \right) $ for $k = 2, 3, 4$.
\end{enumerate}
which defines $W_1$ as the dominating term in the estimate and thus completing the proof. Since the terms $W_2$ and $W_3$ are similar, we shall illustrate the properties of $W_2$ and those of $W_3$ can be derived similarly. 

The term $W_1$ can be expressed as
\begin{eqnarray*}
 W_1 &=& \dmsum_{a,b = -M}^M \xi_n(a,b) \frac{1}{n_{a,b}} \dmsum_{t,s \in \mathcal{A}}  \left(Y_{t+a} - \overline{Y}_n^* \right)' Y_s \left(Y_{s+b} - \overline{Y}_n^* \right)' Y_t \\
 &=& \dmsum_{a,b = -M}^M \xi_n(a,b)\frac{1}{n_{a,b}} \dmsum_{t,s \in \mathcal{A}} \bigg[ Y_{t+a}' Y_s Y_{s+b}' Y_t - \overline{Y}_n^{*'} Y_s Y_{s+b} Y_t \\
 && - Y_{t+a}' Y_s \overline{Y}_n^{*'} Y_t  + \overline{Y}_n^{*'} Y_s \overline{Y}_n^{*'} Y_t \bigg] \\
\end{eqnarray*}
Define 
\begin{eqnarray*}
 W_{11} := \dmsum_{a,b = -M}^M \xi_n(a,b)\frac{1}{n_{a,b}} \dmsum_{t,s \in \mathcal{A}} Y_{t+a}' Y_s Y_{s+b}' Y_t,
\end{eqnarray*}
and $W_{12}, W_{13}$ and $W_{14}$ are defined similarly with the summand in the above expression replaced by $\overline{Y}_n^{*'} Y_s Y_{s+b} Y_t$, $ Y_{t+a}' Y_s \overline{Y}_n^{*'} Y_t$ and $\overline{Y}_n^{*'} Y_s \overline{Y}_n^{*'} Y_t$ respectively. The terms $W_{1k}, k = 1, 2, 3, 4$, when studied individually, yield ratio consistency of the estimate.

\textbf{Properties of $W_{11}$:}

The term $W_{11}$, owing to the construction of the indexing set $\mathcal{A}(a,b)$, is unbiased for $\varmn$,
\begin{align}
 E \left(W_{11} \right) &= E \left[ \dmsum_{a,b = -M}^M \xi_n(a,b)\frac{1}{n_{a,b}} \dmsum_{t,s \in \mathcal{A}} Y_{t+a}' Y_s Y_{s+b}' Y_t \right] \nonumber \\
 &= \dmsum_{a,b = -M}^M \xi_n(a,b)\frac{1}{n_{a,b}} \dmsum_{t,s \in \mathcal{A}} E \left[Y_{t+a}' Y_s Y_{s+b}' Y_t \right] \nonumber\\
 &=  \dmsum_{a,b = -M}^M \frac{\xi_n(a,b)}{n_{a,b}} \msum_{t,s \in \mathcal{A}} \left( \otrab{a}{b} + \otrab{t+a-s-b}{t-s} \right) \nonumber \\
 &= \varmn + \dmsum_{a,b = -M}^M \xi_n(a,b)\frac{1}{n_{a,b}}\otrab{t+a-s-b}{t-s}  \nonumber \\
 &=  \varmn, \label{eqn:ew11}
\end{align}
where the last equality is due to the choice of elements in the indexing set $\mathcal{A}(a,b)$.

Applying Cauchy-Schwarz inequality, the variance of $W_{11}$ can be bounded above as
\begin{eqnarray}
 Var \left(W_{11} \right) &=& Var \left[ \dmsum_{a,b = -M}^M \xi_n(a,b)\frac{1}{n_{a,b}} \dmsum_{t,s \in \mathcal{A}} Y_{t+a}' Y_s Y_{s+b}' Y_t\right] \nonumber \\
 & \leq & \left[ \dmsum_{a,b = -M}^M \xi_n(a,b) \frac{1}{n_{a,b}} \sqrt{ Var \left(  \dmsum_{t,s \in \mathcal{A}} Y_{t+a}' Y_s Y_{s+b}' Y_t \right)} \right]^2,  \label{eqn:varw11}
\end{eqnarray}
since $\xi_n(a,b) \geq 0$. For any $-M \leq a,b \leq M$, $E \left[ \dmsum_{t,s \in \mathcal{A}} Y_{t+a}' Y_s Y_{s+b}' Y_t \right] = n_{a,b} \otrab{a}{b}$, therefore
\begin{align*}
& Var (  \dmsum_{t,s \in \mathcal{A}} Y_{t+a}' Y_s Y_{s+b}' Y_t ) + n_{a,b}^2 tr^2 \left(\Gamma(a) \Gamma(b) \right) = E[ \dmsum_{t,s \in \mathcal{A}} Y_{t+a}' Y_s Y_{s+b}' Y_t ]^2 \\
&= \dmsum_{(t_1,s_1),(t_2,s_2) \in \mathcal{A}} E \left[ Y_{t_1+a}' Y_{s_1} Y_{s_1+b}' Y_{t_1} Y_{t_2+a}' Y_{s_2} Y_{s_2+b}' Y_{t_2}\right].
\end{align*}
For any two pairs of indices $\left(t_1,s_1\right),\left(t_1,s_1\right) \in \mathcal{A}$, from the construction of the indexing set $\mathcal{A}$, we have
\begin{eqnarray}
E \left[ Y_{t_1+a}' Y_{s_1} Y_{s_1+b}' Y_{t_1} Y_{t_2+a}' Y_{s_2} Y_{s_2+b}' Y_{t_2}\right]  &=&  tr^2 \left(\Gamma(a) \Gamma(b) \right) \nonumber \\
 & + &  \dmsum_{\left(q_1,q_2,q_3,q_4\right)} \otrab{q_1}{q_2} \otrab{q_3}{q_4} \nonumber \\
 & + &   \dmsum_{\left(r_1,r_2,r_3,r_4\right)} tr \left( \Gamma(r_1) \Gamma(r_2) \Gamma(r_3) \Gamma(r_4) \right), \label{eqn:eightprods}
\end{eqnarray}
where the coefficients $\left( q_1,q_2,q_3,q_4 \right)$ belong to the set
\begin{eqnarray*}
\mathcal{C}_1 &=& \bigg\{ \left(t_1-t_2,s_2-s_1,t_1-t_2,s_2-s_1 \right), ( t_1-t_2,s_2-s_1, \\
&& t_1-s_2-b,t_2-s_1-b),\left(t_1-t_2,s_2-s_1,t_1-s_2+a,t_2-s_1+a\right), \\
&& (t_1-s_2+a,t_2-s_1+a,t_1-s_1-b,t_2-s_1-b), (s_2-t_1+b,s_1-t_2+b, \\
&& t_1-t_2-a,s_1-s_2-b), (s_2-t_1+b,s_1-t_2+b,s_1-t_2+b-a,s_2-t_1),\\
&& \left(t_2-s_1,t_1-s_2+a-b,t_1-t_2-a,s_1-s_2-b\right), \\
&& \left(t_2-s_1,t_1-s_2+a-b,t_2-s_1+a-b,t_1-s_2\right)\bigg\}.
\end{eqnarray*}
and $\left( r_1,r_2,r_3,r_4 \right)$ belongs to the set
\begin{eqnarray*}
\mathcal{C}_2 &= &\bigg\{  \left(a,s_1-s_2+b,a,s_2-s_1+b\right),  \left(b,t_1-t_2+a,b,t_1-t_2-a\right) \\
&&  \left(-a,t_1-s_2+a-b,b,s_2-t_1\right),  \left(-a,t_2-s_1+a-b,b,s_1-t_2\right),\\
&&  \left(a,t_1-s_2+a,-a,s_2-s_1\right),  \left(b,t_2-t_1,-b,t_1-t_2\right), \\
&&   \left(a,t_1-s_2+a,b,s_2-t_1+b\right),  \left(a,t_2-s_1+a,b,s_1-t_2+b\right), \\
&&  \left(t_2-s_1+a,s_1-s_2+b,s_2-t_1+b,t_1-t_2+a\right),  (t_2-s_1+a,s_1-s_2+b, \\
&&  t_2-t_1, t_1-s_2+a-b), \left(t_2-s_1+a,t_1-s_2,s_2-s_1,t_1-t_2+a\right), \\
&&  \left(t_2-s_1+a,t_1-s_2,t_2-s_1-b,t_1-s_2+a-b\right), \\
&&  \left(s_2-s_1,s_1-t_2+b-a,s_2-t_1+b,t_1-t_2+a\right),  (s_2-s_1,s_1-t_2+b-a,\\
&&  t_2-t_1,t_1-s_2+a-b),\left(s_2-s_1,t_1-t_2-a,s_2-s_1,t_1-t_2+a\right), \\
&&  \left(s_2-s_1,t_1-t_2-a,t_2-s_1-b,t_1-s_2+a-b\right), (t_2-t_1,s_1-s_2+b, \\
&&  s_2-t_1+b,s_1-t_2),  \left(t_2-t_1,s_1-s_2+b,t_2-t_1,s_1-s_2-b\right), \\
&&  \left(t_2-t_1,t_1-s_2, s_2-s_1,s_1-t_2\right),  \left(t_2-t_1,t_1-s_2,t_2-s_1-b,s_1-s_2-b\right), \\
&&  \left(s_2-t_1-a,s_1-t_2+b-a,s_2-t_1+b,s_1-t_2\right),  (s_2-t_1-a,s_1-t_2+b-a, \\
&&  t_2-t_1,s_1-s_2-b), \left(s_2-t_1-a,t_1-t_2-a,s_2-s_1,s_1-t_2\right), \\
&&  \left(s_2-t_1-a,t_1-t_2-a,t_2-s_1-b,s_1-s_2-b\right) \bigg\}.
\end{eqnarray*}

The above indexing sets are obtained by calculating the expected values of products of quadratic forms using results from Magnus \cite{Magnus}. For any $\left(q_1,q_2,q_3,q_4\right)$ with $-M \leq q_i \leq M$, denote by $n_{q_1q_2q_3q_4}$ the number of elements of the set $\mathcal{C}_1$ which are equal to the specified quadruple of indices. Let $n_{r_1r_2r_3r_4}$ denote the corresponding number for elements in $\mathcal{C}_2$. Then we have
\begin{eqnarray*}
 n_{q_1q_2q_3q_4} = \dmsum_{\left(w_1,w_2,w_3,w_4\right) \in \mathcal{C}_1}\# \left\{  w_i = q_i : 1 \leq i \leq 4 \right\}.
\end{eqnarray*}
Using Lemma \eqref{lem:c1c2}, it can be shown that for every $\left(w_1,w_2,w_3,w_4\right) \in \mathcal{C}_1$, the number of elements which are equal to $\left(q_1,q_2,q_3,q_4\right)$ will be of the order $O \left( n^2\right)$. Since the number of elements in $\mathcal{C}_1$ is finite, we have $n_{q_1q_2q_3q_4} = O \left(n^2\right)$. For elements in $\mathcal{C}_2$, there are quadruples $\left(w_1, w_2, w_3, w_4 \right) \in \mathcal{C}_2$ such that at least one of the $w_i$'s are fixed. Therefore the number of elements will be of order $n^3$ and since $\mathcal{C}_2$ has finitely many quadruples, we have $n_{r_1r_2r_3r_4} = O \left(n^3 \right)$. Therefore the quantity in \eqref{eqn:eightprods} can be simplified to obtain
\begin{eqnarray*}
&& Var\left(  \dmsum_{t,s \in \mathcal{A}} Y_{t+a}' Y_s Y_{s+b}' Y_t \right) = \dmsum_{-M \leq r_1,r_2,r_3,r_4 \leq M} n_{r_1r_2r_3r_4} tr \left( \Gamma(r_1) \Gamma(r_2) \Gamma(r_3) \Gamma(r_4) \right) \\
&& +\dmsum_{-M \leq q_1,q_2,q_3,q_4 \leq M } n_{q_1q_2q_3q_4} \otrab{q_1}{q_2} \otrab{q_3}{q_4} 
 \end{eqnarray*}
 
Before studying the asymptotic properties of the variance of $W_{11}$, we establish the following equivalent expression 
\begin{eqnarray*}
 \frac{Var \left(W_{11} \right)}{\varmn^2} &=& \frac{Var \left(W_{11} \right)}{\left(n^{-2} tr\left(\Omega_n^2 \right)\right)^2} \frac{\left(n^{-2} tr \left(\Omega_n^2 \right)\right)^2}{\varmn^2},
 \end{eqnarray*}
where from Proposition \eqref{prop:simplexhis} we have $\frac{\varmn}{n^{-2} tr\left(\Omega_n^2 \right)} = 2 + o(1)$ and
 \begin{eqnarray*}
 \frac{Var \left(W_{11} \right)}{n^{-4} tr^2\left(\Omega_n^2 \right)}& \leq & \left[ \dmsum_{a,b = -M}^M \frac{\xi_n(a,b)}{n^{-2}tr \left(\Omega_n^2\right)}  \sqrt{ \frac{1}{n_{a,b}^2} Var \left(  \dmsum_{t,s \in \mathcal{A}} Y_{t+a}' Y_s Y_{s+b}' Y_t \right)} \right]^2 \\
 &=&  \left[ \dmsum_{a,b = -M}^M n^2 \xi_n(a,b)\sqrt{ \frac{1}{n_{a,b}^2 tr^2 \left(\Omega_n^2 \right)} Var \left(  \dmsum_{t,s \in \mathcal{A}} Y_{t+a}' Y_s Y_{s+b}' Y_t \right)} \right]^2.
\end{eqnarray*}
For any two sequences $\left\{a_n\right\}$ and $\left\{b_n\right\}$, define $a_n \sim b_n$ if $\lim \frac{a_n}{b_n} = 1$. For any $-M \leq a,b \leq M$, $n^2 \xi_n(a,b) \sim 1$ and
\begin{eqnarray}
 &&\frac{1}{n_{a,b}^2 tr^2 \left(\Omega_n^2 \right)} Var \left(  \dmsum_{t,s \in \mathcal{A}} Y_{t+a}' Y_s Y_{s+b}' Y_t \right) \nonumber  \\
 &&= \frac{1}{n_{a,b}^2}  \dmsum_{q_1,q_2,q_3,q_4 = -M}^M   n_{q_1q_2q_3q_4} \frac{\otrab{q_1}{q_2} \otrab{q_3}{q_4}}{tr^2 \left( \Omega_n^2 \right)} \nonumber \\
 && + \frac{1}{n_{a,b}^2} \dmsum_{r_1,r_2,r_3,r_4 = -M}^M  n_{r_1r_2r_3r_4}  \frac{tr \left( \Gamma(r_1) \Gamma(r_2) \Gamma(r_3) \Gamma(r_4) \right)}{tr^2 \left( \Omega_n^2 \right)} \nonumber \\
 &&= \frac{1}{n_{a,b}^2} \dmsum_{q_1,q_2,q_3,q_4= -M}^M  n_{q_1q_2q_3q_4} \frac{\otrab{q_1}{q_2}}{ tr\left( \Omega_n^2 \right)}  \frac{\otrab{q_3}{q_4}}{ tr \left( \Omega_n^2 \right)} \nonumber \\
 && +  \frac{1}{n_{a,b}^2} \dmsum_{r_1,r_2,r_3,r_4 = -M}^M  n_{r_1r_2r_3r_4}  \frac{tr \left( \Gamma(r_1) \Gamma(r_2) \Gamma(r_3) \Gamma(r_4) \right)}{tr^2 \left( \Omega_n^2 \right)}. \label{eqn:varratio1}
\end{eqnarray}
Since $\frac{n_{a,b}}{n^2} \rightarrow 1$, $p = O(n)$, $n_{q_1q_2q_3q_4} = O \left(n^2 \right)$, from Lemma \eqref{lem:gammaproplemma} and \eqref{eqn:condition1}, $n_{r_1r_2r_3r_4} = O \left(n^3 \right)$ and $M = O \left(n^{\frac{1}{8}} \right)$, we have 
\begin{eqnarray*}
 \frac{1}{n_{a,b}^2} \dmsum_{q_1,q_2,q_3,q_4= -M}^M  n_{q_1q_2q_3q_4} \frac{\otrab{q_1}{q_2}}{ tr\left( \Omega_n^2 \right)}  \frac{\otrab{q_3}{q_4}}{ tr \left( \Omega_n^2 \right)}  & =& \frac{1}{n^4} M^4 O \left(n^2 \right) o\left(n^{\frac{1}{2}} \right) \left(n^{\frac{1}{2}} \right), \\
  \frac{1}{n_{a,b}^2} \dmsum_{r_1,r_2,r_3,r_4 = -M}^M  n_{r_1r_2r_3r_4}  \frac{tr \left( \Gamma(r_1) \Gamma(r_2) \Gamma(r_3) \Gamma(r_4) \right)}{tr^2 \left( \Omega_n^2 \right)} & = & \frac{1}{n^4} M^4 O \left(n^3 \right) o(1),
\end{eqnarray*}
where the quantities on the right hand side are equal to $\frac{1}{n^4} O \left(n^{\frac{4}{8}} \right) o \left( n^3 \right) = o \left( n^{-\frac{1}{2}} \right)$ and $\frac{1}{n^4} O \left(n^{\frac{4}{8}} \right) O \left(n^3 \right) o (1) = o \left(n^{-\frac{1}{2}} \right)$ respectively. Therefore, we have 
\begin{eqnarray}
 \frac{1}{n_{a,b}^2 tr^2 \left(\Omega_n^2 \right)} Var \left(  \dmsum_{t,s \in \mathcal{A}} Y_{t+a}' Y_s Y_{s+b}' Y_t \right) = o \left(n^{-\frac{1}{2}} \right),
 \label{eqn:varw11rate}
\end{eqnarray}
which gives $\frac{Var \left(W_{11} \right)}{n^{-4} tr^2 \left(\Omega_n^2 \right)} = o \left(n^{-\frac{1}{2}} \right) M^4 = o \left(n^{-\frac{1}{2}} \right) O \left(n^{\frac{1}{2}} \right) = o(1)$. This result combined with \eqref{eqn:ew11} implies $W_{11}$ is consistent for $\varmn$, i.e
\begin{eqnarray}
 \frac{W_{11}}{\varmn} \stackrel{p}{\longrightarrow} 1.
 \label{eqn:w11consistent}
\end{eqnarray}

\noindent \textbf{Properties of $W_{12}$}

Looking at the expected value of $W_{12}$,
\begin{eqnarray*}
 E \left[ W_{12} \right] &=& E \left[\dmsum_{a,b = -M}^M \xi_n(a,b) \frac{1}{n_{a,b}} \dmsum_{t,s \in \mathcal{A}} \overline{Y}_n^{*'} Y_s Y_{s+b} Y_t\right], \\
 &=& \dmsum_{a,b = -M}^M \xi_n(a,b) \frac{1}{n_{a,b}} \dmsum_{t,s \in \mathcal{A}} \frac{1}{m_{\mathcal{B}}} \msum_{w \in \mathcal{B}} E \left[Y_w' Y_s Y_{s+b}' Y_t\right],
\end{eqnarray*}
where $m_{\mathcal{B}}$ is the cardinality of the indexing set $\mathcal{B}$ and $E \left[Y_w' Y_s Y_{s+b}' Y_t\right] = \otrab{w-t}{b} + \otrab{w-s-b}{t-s}$, which is equal to zero by construction of $\mathcal{B}(t,s)$. The variance of $W_{12}$ can be bounded above as
\begin{eqnarray}
 Var \left(W_{12}\right) &=& Var \left[\dmsum_{a,b = -M}^M \xi_n(a,b) \frac{1}{n_{a,b}} \dmsum_{t,s \in \mathcal{A}} \overline{Y}_n^{*'} Y_s Y_{s+b} Y_t\right] \nonumber \\
 & \leq & \left[ \dmsum_{a,b = -M}^M \xi_n(a,b) \frac{1}{n_{a,b}} \sqrt{ Var \left( \dmsum_{t,s \in \mathcal{A}} \overline{Y}_n^{*'} Y_s Y_{s+b} Y_t \right)} \right]^2 \nonumber \\
 & \leq & \left[ \dmsum_{a,b = -M}^M \xi_n(a,b) \frac{1}{n_{a,b}} \sqrt{ \left[\dmsum_{t,s \in \mathcal{A}} \sqrt{Var \left(  \overline{Y}_n^{*'} Y_s Y_{s+b} Y_t\right)} \right]^2} \right]^2 \nonumber \\
 &=& \left[ \dmsum_{a,b = -M}^M \xi_n(a,b) \frac{1}{n_{a,b}} \dmsum_{t,s \in \mathcal{A}} \sqrt{Var \left(  \overline{Y}_n^{*'} Y_s Y_{s+b} Y_t\right)}  \right]^2, \label{eqn:w12var}
\end{eqnarray}
where the two inequalities are from $Var \left(\mathop{\sum} a_i X_i \right) \leq \left( \mathop{\sum} a_i \sqrt{Var(X_i)} \right)^2$ for $a_i > 0$. For any $t,s \in \mathcal{A}$,
\begin{eqnarray*}
 Var \left(  \overline{Y}_n^{*'} Y_s Y_{s+b}' Y_t\right) &=& E \left( \overline{Y}_n^{*'} Y_s Y_{s+b}' Y_t \right)^2 \\
 &=& \frac{1}{m_{\mathcal{B}}^2} \left[ tr \left( \Omega_{n \mathcal{B}} \Gamma(0) \right) tr\left(\Gamma(0)^2 \right) + 2 \, tr \left[ \Omega_{n \mathcal{B}} \Gamma(b) \Gamma(0) \Gamma(-b) \right]\right] \\
 &=& \frac{1}{m_{\mathcal{B}}^2} \msum_{h=-M}^M  a_{n, \mathcal{B}}(h) \bigg[  tr \left( \Gamma(h) \Gamma(0) \right) tr\left(\Gamma(0)^2 \right) \\
 && + 2 \, tr \left[ \Gamma(h) \Gamma(b) \Gamma(0) \Gamma(-b) \right]  \bigg],
\end{eqnarray*}
where $m_{\mathcal{B}} \sim  n$ is the cardinality of $\mathcal{B}$ and $a_{n,\mathcal{B}}(h) = 1 + o(1), -M \leq h \leq M$ from Lemma \eqref{lem:btslemma}. The coefficients $a_{n,\mathcal{B}}(h)$ are similar to $b_i$'s defined for $\Omega_n$. Therefore from \eqref{eqn:w12var} and Proposition \eqref{prop:simplexhis}, we have
\begin{eqnarray*}
 \frac{Var \left(W_{12}\right)}{n^{-4}tr^2 \left(\Omega_n^2 \right)} &\leq &  \left[ \dmsum_{a,b = -M}^M n^2 \xi_n(a,b) \frac{1}{n_{a,b}} \dmsum_{t,s \in \mathcal{A}} \sqrt{ \frac{Var \left(  \overline{Y}_n^{*'} Y_s Y_{s+b} Y_t\right)}{tr^2 \left(\Omega_n^2 \right)}}  \right]^2,  \\
\end{eqnarray*}
and for any $(t,s) \in \mathcal{A}$,
\begin{eqnarray*}
  \frac{Var \left(  \overline{Y}_n^{*'} Y_s Y_{s+b} Y_t\right)}{tr^2 \left(\Omega_n^2 \right)}
  &=& \frac{1}{tr^2 \left(\Omega_n^2 \right)} \frac{1}{m_{\mathcal{B}}^2} \msum_{h=-M}^M  a_{n, \mathcal{B}}(h) \bigg[  tr \left( \Gamma(h) \Gamma(0) \right) tr\left(\Gamma(0)^2 \right) \\
 && \, + 2 \, tr \left[ \Gamma(h) \Gamma(b) \Gamma(0) \Gamma(-b) \right]  \bigg], \\
 &=& \frac{n}{m_{\mathcal{B}}^2} \msum_{h=-M}^M a_{n, \mathcal{B}}(h) \bigg[  \frac{tr \left( \Gamma(h) \Gamma(0) \right)}{\sqrt{n} tr \left(\Omega_n^2 \right)} \frac{tr\left(\Gamma(0)^2 \right)}{ \sqrt{n} tr \left(\Omega_n^2 \right)}\\
 && \, + 2 \, \frac{tr \left[ \Gamma(h) \Gamma(b) \Gamma(0) \Gamma(-b) \right]}{n tr^2 \left(\Omega_n^2 \right)}  \bigg] \\
 &=& O \left(\frac{1}{n} \right) (1 + o(1)) \left(2M+1 \right) \left(o(1) 1 + o \left(n^{-1} \right) \right),
\end{eqnarray*}
since $M = O \left(n^{\frac{1}{8}} \right)$ and the summands in the above expression are $o \left(n^{-\frac{7}{8}} \right)$ according to Lemma \eqref{lem:gammaproplemma} and the assumption in \eqref{eqn:condition1}. Therefore, it can be established that $\frac{Var \left(W_{12}\right)}{n^{-4}tr^2 \left(\Omega_n^2 \right)} = o(n^{-\frac{3}{8}})$ and by Proposition \eqref{prop:simplexhis}, we can conclude that
\begin{eqnarray*}
 \frac{Var \left(W_{12}\right)}{\varmn^2} \longrightarrow 0 \, \, \mbox{ as } n \rightarrow \infty.
\end{eqnarray*}

By similar arguments we can establish convergence of the remaining terms, $W_{13} = o_p\left(\varmn\right)$ and $W_{14} = o_p\left(\varmn\right)$. Combining the properties of $W_{1k}, k = 1, 2, 3 ,4$, $W_1$ is ratio consistent for $\varmn$. 

Since $W_2$ and $W_3$ comprise of product of two quadratic terms where one of the terms is equal to $\mu$, the calculations for consistency follow along the same lines. By the same procedure used to establish consistency of $W_1$, we can express $W_2$ as
\begin{eqnarray*}
 W_2 &=& \dmsum_{a,b = -M}^M \xi_n(a,b) \frac{1}{n_{a,b}} \dmsum_{t,s \in \mathcal{A}} \left( Y_{t+a} - \overline{Y}_n^* \right)' \mu \left( Y_{s+b} - \overline{Y}_n^* \right)' Y_t \\
 &=& \dmsum_{a,b = -M}^M \xi_n(a,b) \frac{1}{n_{a,b}}  \dmsum_{t,s \in \mathcal{A}} \bigg(  Y_{t+a}' \mu Y_{s+b} Y_t - \mu' \overline{Y}_n^*  Y_{s+b}' Y_t \\
 & & - Y_{t+a}' \mu Y_t' \overline{Y}_n^*  + \mu' \overline{Y}_n^*  Y_t' \overline{Y}_n^* \bigg) 
\end{eqnarray*}
Define $W_{21} = \dmsum_{a,b = -M}^M \xi_n(a,b) \frac{1}{n_{a,b}}  \dmsum_{t,s \in \mathcal{A}}  Y_{t+a}' \mu Y_{s+b} Y_t$ and $W_{22}, W_{23}$ and $W_{24}$ are defined similarly with the summands replaced by $\mu' \overline{Y}_n^*  Y_{s+b}' Y_t, Y_{t+a}' \mu Y_t' \overline{Y}_n^*$ and $\mu' \overline{Y}_n^*  Y_t' \overline{Y}_n^*$ respectively.

\noindent \textbf{Properties of $W_{21}$}:

The expected value of $W_{21}$ is equal to zero because it comprises a linear combination of third-order moments of zero mean Gaussian distrbution. The variance of $W_{21}$ can be calculated as
\begin{eqnarray*}
 Var \left(W_{21} \right) &=& Var \left[ \dmsum_{a,b = -M}^M \xi_n(a,b) \frac{1}{n_{a,b}}  \dmsum_{t,s \in \mathcal{A}}  Y_{t+a}' \mu Y_{s+b} Y_t \right] \\
 & \leq & \left[ \dmsum_{a,b = -M}^M \xi_n(a,b) \frac{1}{n_{a,b}}\sqrt{ Var \left(\dmsum_{t,s \in \mathcal{A}}  Y_{t+a}' \mu Y_{s+b} Y_t \right)}\right]^2  \\
 & = & \left[ \dmsum_{a,b = -M}^M \xi_n(a,b) \frac{1}{n_{a,b}} \left( \dmsum_{t,s \in \mathcal{A}}  \sqrt{ Var \left( Y_{t+a}' \mu Y_{s+b} Y_t \right)} \right)\right]^2.
\end{eqnarray*}
For any $(t,s) \in \mathcal{A}$, the variance of $Y_{t+a}' \mu Y_{s+b} Y_t$ is
\begin{eqnarray*}
 Var \left( Y_{t+a}' \mu Y_{s+b} Y_t \right) &=& E  \left[ Y_{t+a}' \mu Y_{s+b} Y_t \right]^2 \\
 &=& \frac{1}{2} \left(\mu' \Gamma(0) \mu \right) tr \left(\Gamma(0)^2 \right)  + \frac{1}{4} \mu' \Gamma(-a) \Gamma(0) \Gamma(a) \mu
\end{eqnarray*}
which implies
\begin{eqnarray*}
 && \frac{Var \left(W_{21} \right)}{n^{-4}tr^2 \left( \Omega_n^2 \right)} \leq  \left[ \dmsum_{a,b = -M}^M \frac{\xi_n(a,b)}{n^{-2} tr \left( \Omega_n^2 \right)} \frac{1}{n_{a,b}} \left( \dmsum_{t,s \in \mathcal{A}}  \sqrt{ Var \left( Y_{t+a}' \mu Y_{s+b} Y_t \right)} \right)\right]^2 \\
 && = \left[ \dmsum_{a,b = -M}^M \frac{\xi_n(a,b)}{n^{-2} tr \left( \Omega_n^2 \right)} \sqrt{  \frac{1}{2} \left(\mu' \Gamma(0) \mu \right) tr \left(\Gamma(0)^2 \right)  + \frac{1}{4} \mu' \Gamma(-a) \Gamma(0) \Gamma(a) \mu }\right]^2 \\
 &&= \left[ \dmsum_{a,b = -M}^M \frac{n^2 \xi_n(a,b)}{2} \sqrt{ \frac{2}{\sqrt{n}} \frac{\mu' \Gamma(0) \mu}{n^{-1}tr \left( \Omega_n^2 \right)}  \frac{tr \left(\Gamma(0)^2 \right)}{\sqrt{n} tr \left( \Omega_n^2 \right)}  + \frac{\mu' \Gamma(-a) \Gamma(0) \Gamma(a) \mu}{tr^2 \left( \Omega_n^2 \right)}} \right]^2.
\end{eqnarray*}
The number of terms in the above summation is $O \left(n^{\frac{1}{4}} \right)$ while the summands can be shown to be $o \left(n^{-\frac{1}{4}}\right)$ by \eqref{eqn:alternativemean}, Lemma \eqref{lem:gammaproplemma} and Lemma \eqref{lem:meancond1}. Therefore, the above quantity converges to zero which implies $W_{21} = o_p \left( \varmn \right)$. By the same argument, we can establish $W_{2k} = o_p\left( \varmn \right)$ for $k = 2, 3$ and $4$, which implies $W_2 = o_p \left( \varmn \right)$. Since $W_3$ involves terms similar to $W_3$, we can repeat the argument given above to conclude that $W_3 = o_p \left( \varmn \right)$.

To show $W_4 = o_p \left( \varmn \right)$, we need the following additional property which results from the local alternative condition mentioned in \eqref{eqn:alternativemean}. Expressing $W_4$ as a sum of four terms as before, we have
\begin{eqnarray*}
W_4 & = &  \dmsum_{a,b=-M}^M \xi_n(a,b) \frac{1}{n_{a,b}} \dmsum_{t,s \in \mathcal{A}} \left( Y_{t+a} - \overline{Y}_n^* \right)' \mu \left( Y_{s+b} - \overline{Y}_n^* \right)' \mu \\
&=& \dmsum_{a,b=-M}^M \xi_n(a,b) \frac{1}{n_{a,b}} \dmsum_{t,s \in \mathcal{A}} \bigg( Y_{t+a}' \mu \mu' Y_{s+b} - \mu'  \overline{Y}_n^*  Y_{s+b}' \mu \\
&& - Y_{t+a}' \mu \mu'\overline{Y}_n^*  + \mu'\overline{Y}_n^* \mu'  \overline{Y}_n^*  \bigg).
\end{eqnarray*}
Define $W_{41} = \dmsum_{a,b = -M}^M \xi_n(a,b) \frac{1}{n_{a,b}}  \dmsum_{t,s \in \mathcal{A}} Y_{t+a}' \mu \mu' Y_{s+b}$ and $W_{42}, W_{43}$ and $W_{44}$ similarly with the summands replaced by $\mu'  \overline{Y}_n^*  Y_{s+b}' \mu,Y_{t+a}' \mu \mu'\overline{Y}_n^*$ and $\mu'\overline{Y}_n^* \mu'  \overline{Y}_n^*$  respectively.

\noindent \textbf{Properties of $W_{41}$}:

By construction of the indexing sets $\mathcal{A}$ and $\mathcal{B}_{t,s}$, $\left( Y_{t+a}, Y_{s+b}, \overline{Y}_n^* \right)$ are mutually independent for any $(t,s) \in \mathcal{A}$. Therefore the first three terms will have zero expected value, $E \left[ Y_{t+a}' \mu \mu' Y_{s+b} \right] = E \left[  \mu'  \overline{Y}_n^*  Y_{s+b}' \mu \right] = E \left[ Y_{t+a}' \mu \mu'\overline{Y}_n^* \right] = 0$ which implies $E \left[ W_{4k} \right] = 0$ for $k = 1, 2, 3$. When $k = 4$, $E \left[\mu'\overline{Y}_n^* \mu'  \overline{Y}_n^* \right] = \mu' Var \left( \overline{Y}_n^* \right) \mu = \frac{1}{m_\mathcal{B}} \mu' \Omega_{n \mathcal{B}}^* \mu$ where $\Omega_{n \mathcal{B}}^*$ is as defined in Lemma \eqref{lem:btslemma}. The quadratic product $\mu' \Omega_{n \mathcal{B}}^* \mu$ is a linear combination of $ \mu' \Gamma(h) \mu $ for $-M \leq h \leq M$ with coefficients given by $\frac{1}{m_{\mathcal{B}}} a_{n, \mathcal{B}}(h) = \frac{1}{n} (1 + o(1))$ from Lemma \eqref{lem:btslemma}. Therefore by lemma \eqref{lem:alternative2} we have $\mu' \Omega_{n \mathcal{B}}^* \mu = o \left(n^{-7/8} tr \left(\Omega_n^2 \right) \right)$ 
which implies $W_{44} =  o_p \left(n^{-1} tr \left(\Omega_n^2 \right) \right)$.

The variance of $W_{4k},k = 1, 2, 3, 4$ can be bounded above as follows
\begin{displaymath}
 Var \left(W_{4k} \right) \leq \left\{
     \begin{array}{ll}
       \left[ \dmsum_{a,b=-M}^M \xi_n(a,b) \left( \mu' \Gamma(0) \mu \right)\right]^2, & k = 1 \\
       \left[ \dmsum_{a,b=-M}^M \xi_n(a,b) \frac{1}{n_{a,b}} \dmsum_{t,s \in \mathcal{A}} \sqrt{\mu' \Gamma(0) \mu \mu' \Omega_{n \mathcal{B}}^* \mu} \right]^2, &  k = 2,3 \\
       \left[ \dmsum_{a,b=-M}^M \xi_n(a,b) \frac{1}{n_{a,b}} \dmsum_{t,s \in \mathcal{A}} \left( \mu' \Omega_{n \mathcal{B}}^*  \mu \right)\right]^2. &  k = 4
     \end{array}
   \right.
\end{displaymath}
Recall from \eqref{eqn:traceest}, the indexing set $\mathcal{A}$ is defined as $\mathcal{A}(a,b) = \{ (t,s) : |t+a-s-b| > M, |t-s| > M \}$, with $n_{a,b} = |\mathcal{A}| \sim n$ and $\xi_n(a,b) \sim \frac{1}{n^2}$.
By Lemma \eqref{lem:alternative2}, the above upper bounds can be shown to be $o_p \left( n^{-4} tr^2 \left( \Omega_n^2 \right) \right)$ which proves $\frac{W_{4k}}{\varmn} \rightarrow 0$ in probability for $k = 1, 2, 3, 4$. This establishes $W_4 = o_p \left(\varmn \right)$, thus completing the proof. 

\qed

\noindent \textbf{Proof of Theorem \eqref{thm:twosample}}

From the expression $\mtxt{M} = \twofam{\mtxt{M}}{1} + \twofam{\mtxt{M}}{2} - 2 \twofam{\mtxt{M}}{3} $, the first two terms come entirely from the first and second sample respectively. Therefore by  Theorem \eqref{thm:clt}, we can establish asymptotic normality of the first two terms. When the third term $\twofam{\mtxt{M}}{3} = \overline{X}^{(1) \prime}_n \overline{X}^{(2)}_n = \mathop{\sum}_{i=1}^{n_1}  \mathop{\sum}_{j=1}^{n_2} X_i^{(1) \prime}  X_i^{(2)} $ is written as the sum of a rectangular array, we can apply the method of triangular arrays as mentioned in the proof of Theorem \eqref{thm:clt} to derive asymptotic normality. Therefore, we have
\begin{eqnarray}
 \frac{\twofam{\mtxt{M}}{i}}{\sqrt{Var \left( \twofam{\mtxt{M}}{i} \right) }} \xrightarrow{\mathcal{D}} \mathcal{N}(0,1) \mbox{ for } i = 1, 2, 3.
\end{eqnarray}
and since the terms are independent, $\mtxt{M}$ is a linear combination of the three terms which establishes asymptotic normality of $\mtxt{M}$.

Consistency of the variance estimator can be proved by extending results from the one sample case. Since $\twofam{\mtxt{M}}{i}, i = 1, 2, 3$ are independent, the variance can be expressed as the sum of variances as $\varmn = Var \left( \twofam{\mtxt{M}}{1} \right) +  Var \left( \twofam{\mtxt{M}}{2} \right) + 4  Var \left( \twofam{\mtxt{M}}{3} \right)$. Ratio consistency of the first two terms can be established using Theorem \eqref{thm:variance}. The third term of the variance is
\begin{eqnarray*}
 Var \left( \twofam{\mtxt{M}}{3} \right) &=& \frac{2}{n_1 n_2} \widehat{tr \left(\twofam{\Omega_n}{1} \twofam{\Omega_n}{2} \right)} \\
 &=& \frac{2}{n_1 n_2} \dmsum_{a,b=-M}^M \left(1 - \frac{|a|}{n_1} \right)\left(1 - \frac{|b|}{n_2} \right) \widehat{tr \left(\twofam{\Gamma}{1}(a) \twofam{\Gamma}{2}(b) \right)},
\end{eqnarray*}
where $\widehat{tr \left(\twofam{\Gamma}{1}(a) \twofam{\Gamma}{2}(b) \right)}$ is constructed as
\begin{eqnarray*} 
&& \dmsum_{a,b = -M}^M  \frac{\xi_n(a,b)}{n_{a,b}} \mathop{\sum \sum}_{\substack{t \in \twofam{\mathcal{A}}{1}(a),\\ s \in \twofam{\mathcal{A}}{2}(b)}} \left(\twofam{X}{1}_{t+a} - \overline{\twofam{X}{1}}_n^* \right)' \twofam{X}{2}_{s} \left(\twofam{X}{2}_{s+b} - \overline{\twofam{X}{2}}_n^* \right)' \twofam{X}{1}_{t} \\
 &=& \dmsum_{a,b = -M}^M \frac{\xi_n(a,b)}{n_{a,b}} \mathop{\sum \sum}_{\substack{t \in \twofam{\mathcal{A}}{1}(a),\\ s \in \twofam{\mathcal{A}}{2}(b)}} \bigg\{ \left(\twofam{Y}{1}_{t+a} - \overline{Y}_n^{(1)*} \right)' \twofam{Y}{2}_s \left(\twofam{Y}{2}_{s+b} - \overline{Y}_n^{(2)*} \right)' \twofam{Y}{1}_t \\
 && +  \left(\twofam{Y}{1}_{t+a} - \overline{Y}_n^{(1)*}\right)' \twofam{\mu}{2} \left(\twofam{Y}{2}_{s+b} - \overline{Y}_n^{(2)*} \right)' \twofam{Y}{1}_{t} \\
 && + \left(\twofam{Y}{1}_{t+a} - \overline{Y}_n^{(1)*} \right)' \twofam{Y}{2}_s \left(\twofam{Y}{2}_{s+b} - \overline{Y}_n^{(2)*} \right)' \twofam{\mu}{1} \\
 && + \left(\twofam{Y}{1}_{t+a} - \overline{Y}_n^{(1)*} \right)' \twofam{\mu}{2} \left(\twofam{Y}{2}_{s+b} - \overline{Y}_n^{(2)*} \right)' \twofam{\mu}{1} \bigg\} \\
 &=& \msum_{k=1}^4 W_k.
\end{eqnarray*}
In the above expression, the observations $\twofam{X}{.}$ are centered to obtain $\twofam{Y}{.}$. Following the calculations in Section \eqref{sec:varianceproof} for establishing the ratio consistency of variance in the one sample case, under the conditions mentioned in Theorem \eqref{thm:twosample}, it is fairly straightforward to prove
\begin{eqnarray}
 \frac{\widehat{Var \left( \twofam{\mtxt{M}}{i} \right)}}{Var \left( \twofam{\mtxt{M}}{i} \right)} \xrightarrow{p} 1 \mbox{  for  } i = 1, 2, 3.
 \label{eqn:var2consistent}
\end{eqnarray}
Combining equation \eqref{eqn:var2consistent} with asymptotic normality of $\mtxt{T}$, by Slutsky's Theorem we have $T_{new} = \frac{\mtxt{M}}{\sqrt{\widehat{Var \left( \mtxt{M} \right)}}}$ to be asymptotically normal for the two sample case.

\qed

\pagebreak
\section{ Supplementary Material: Size and Power Plots}
\label{sec:suppl}
\begin{figure}[htpb]
\centering
 \includegraphics[scale=0.75]{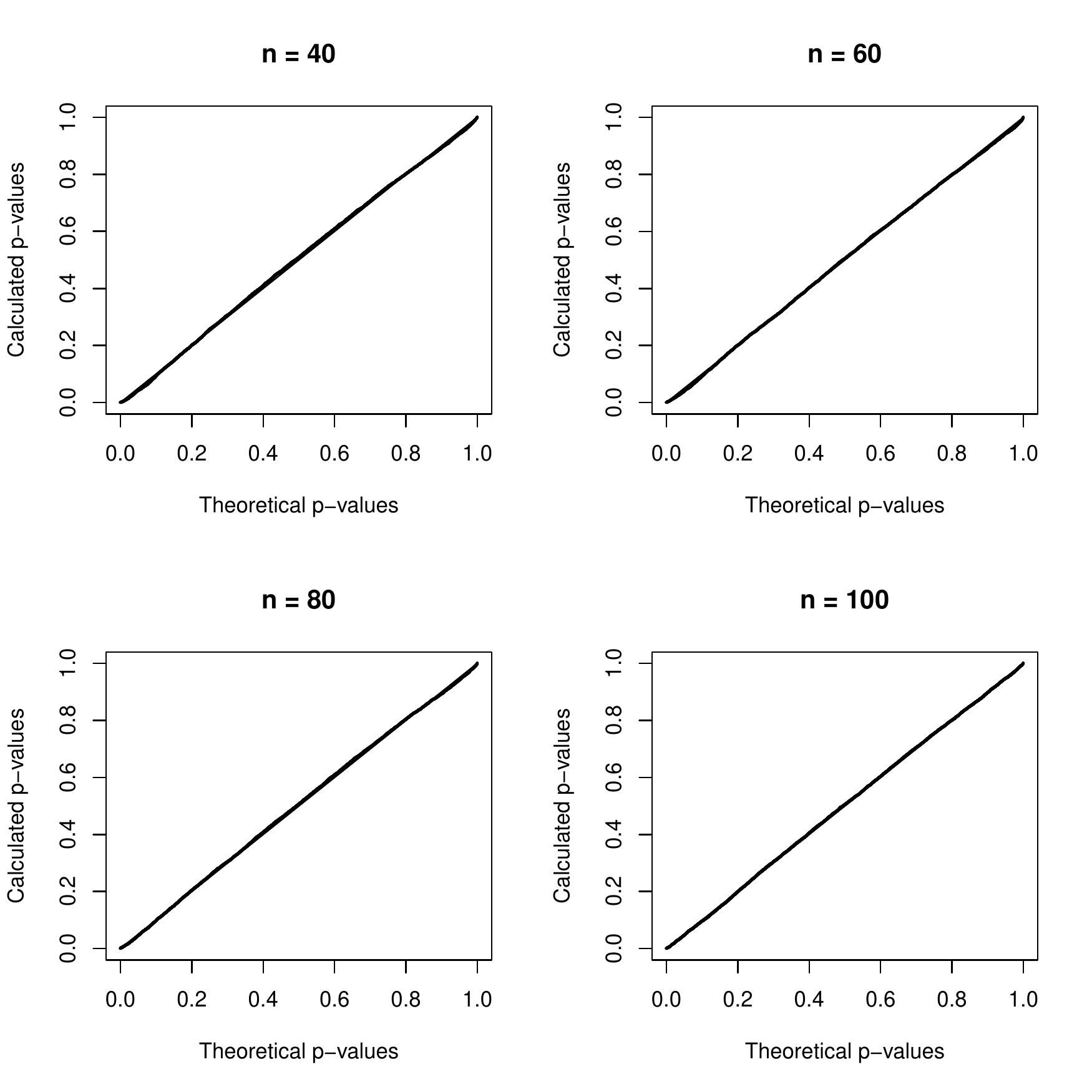}
 \caption{QQ-plot of the p-values calculated from 10,000 randomly generated data sets under model I.}
 \label{fig:qqplotmodel0}
\end{figure}

\begin{figure}[htp]
\centering
 \includegraphics[scale=0.75]{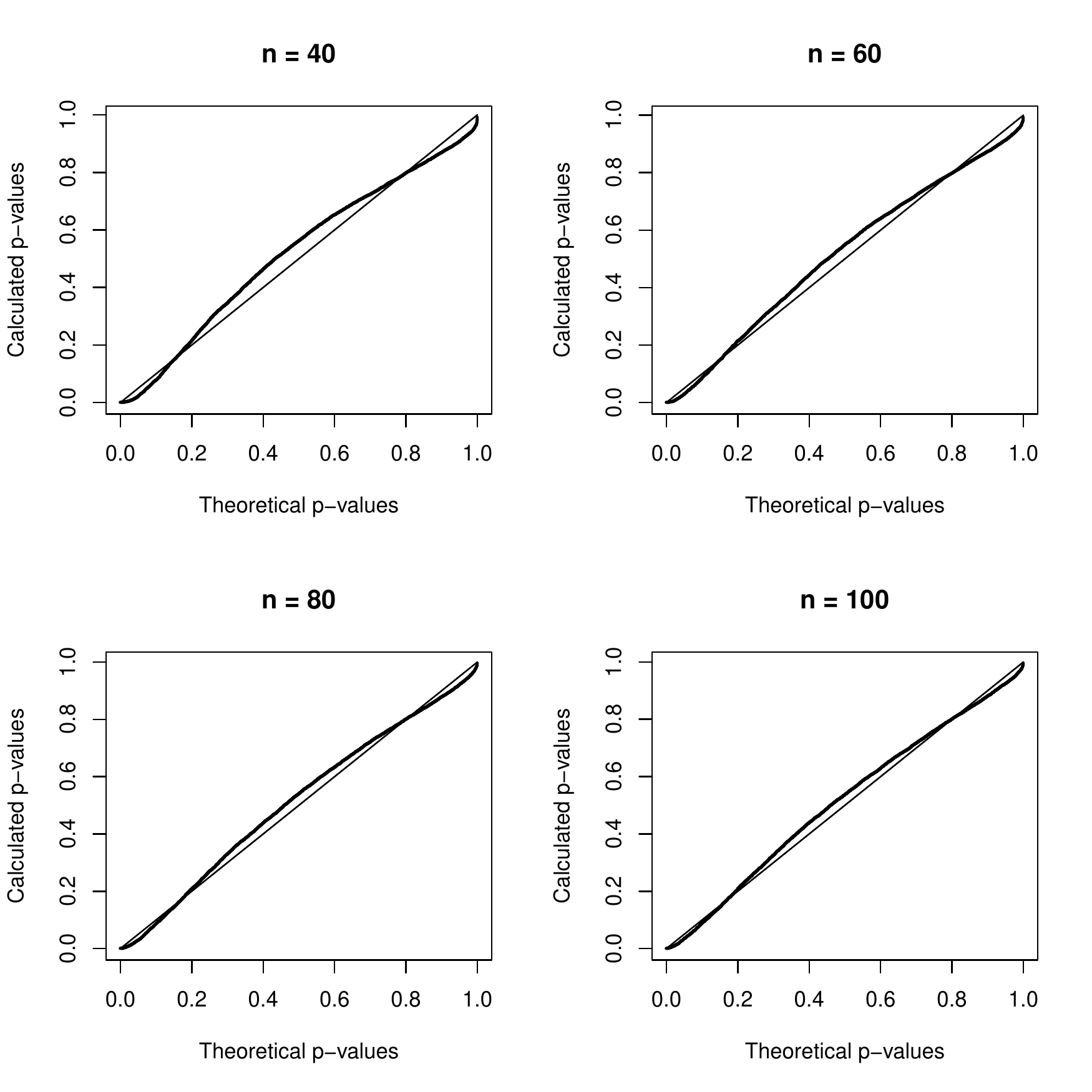}
 \caption{QQ-plot of the p-values calculated from 10,000 randomly generated data sets under model II.}
 \label{fig:qqplotmodel1}
\end{figure}

\begin{figure}[htp]
\centering
 \includegraphics[scale=0.75]{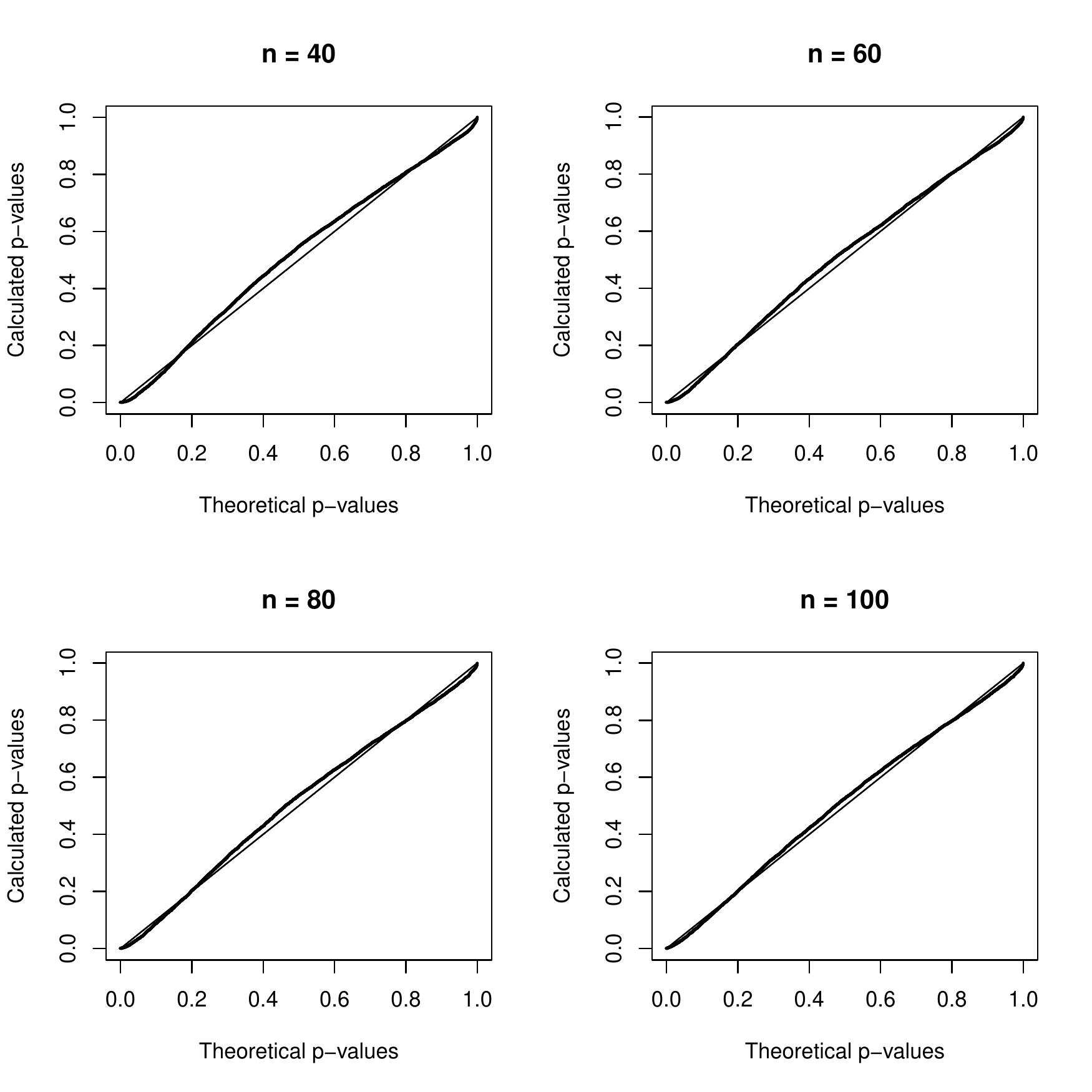}
 \caption{QQ-plot of the p-values calculated from 10,000 randomly generated data sets under model III.}
 \label{fig:qqplotmodel2}
\end{figure}

\begin{figure}[htp]
\centering
 \includegraphics[scale=0.75]{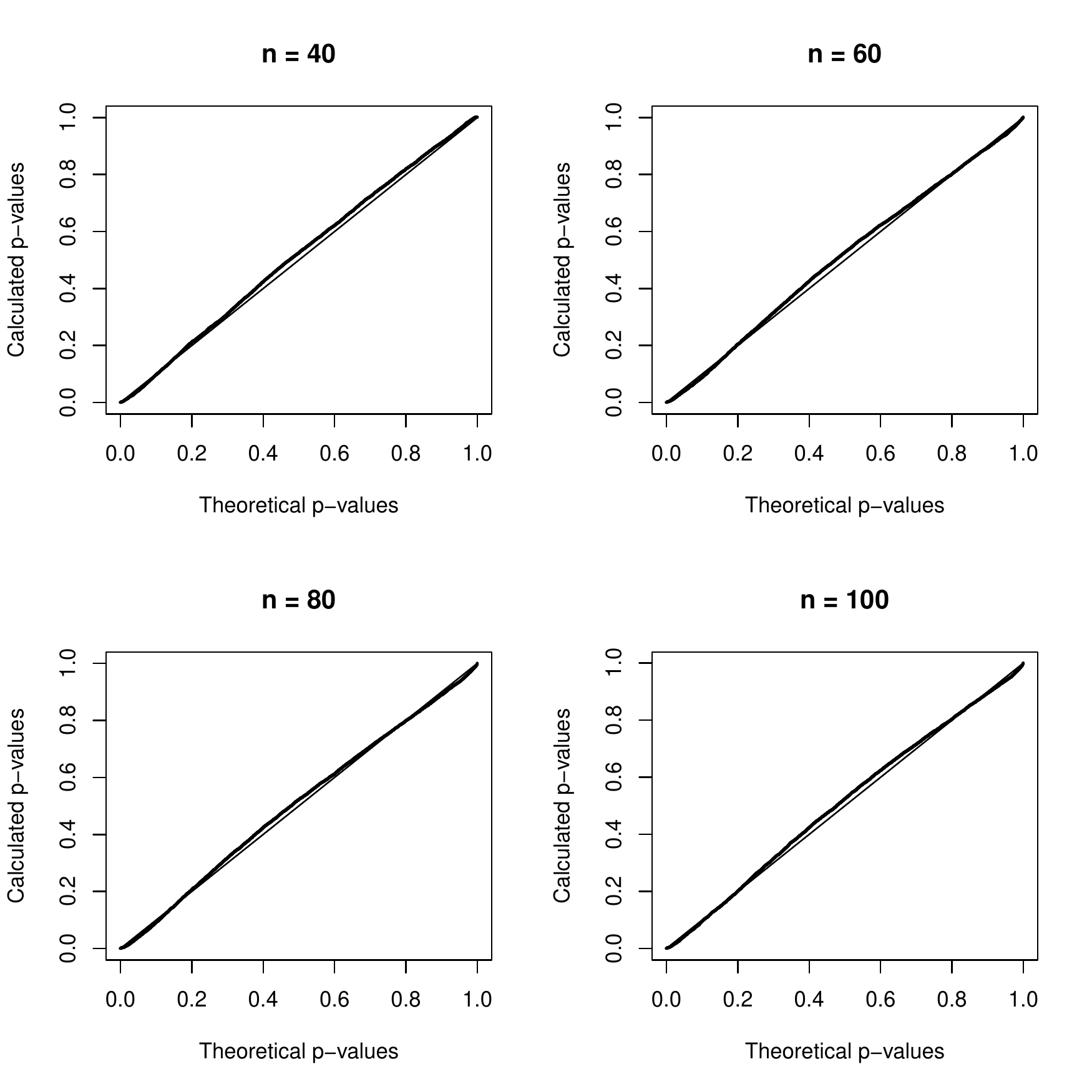}
 \caption{QQ-plot of the p-values calculated from 10,000 randomly generated data sets under model IV.}
 \label{fig:qqplotmode1l3}
\end{figure}

\begin{figure}[htp]
\centering
 \includegraphics[width=0.3\textwidth]{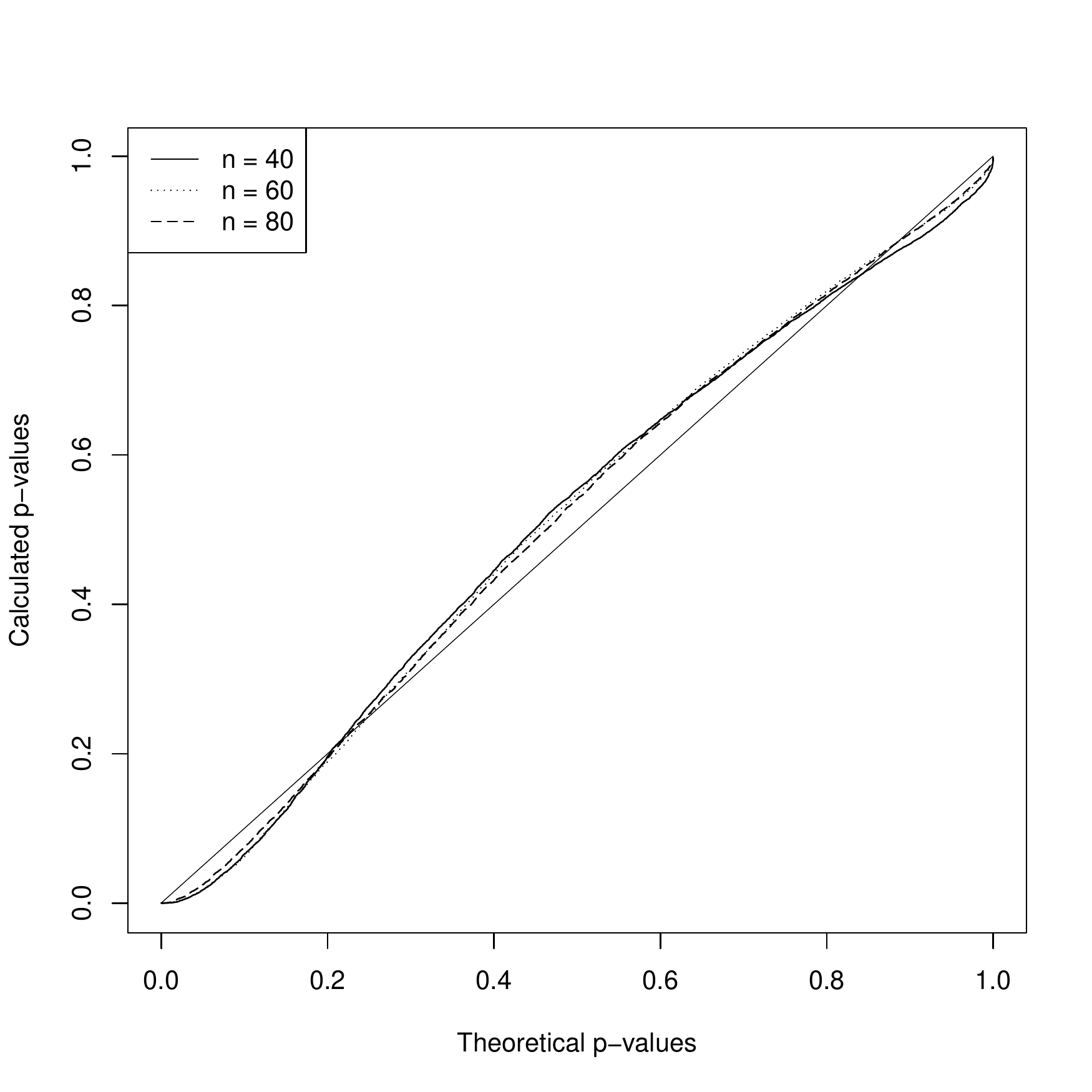} 
 \includegraphics[width=0.3\textwidth]{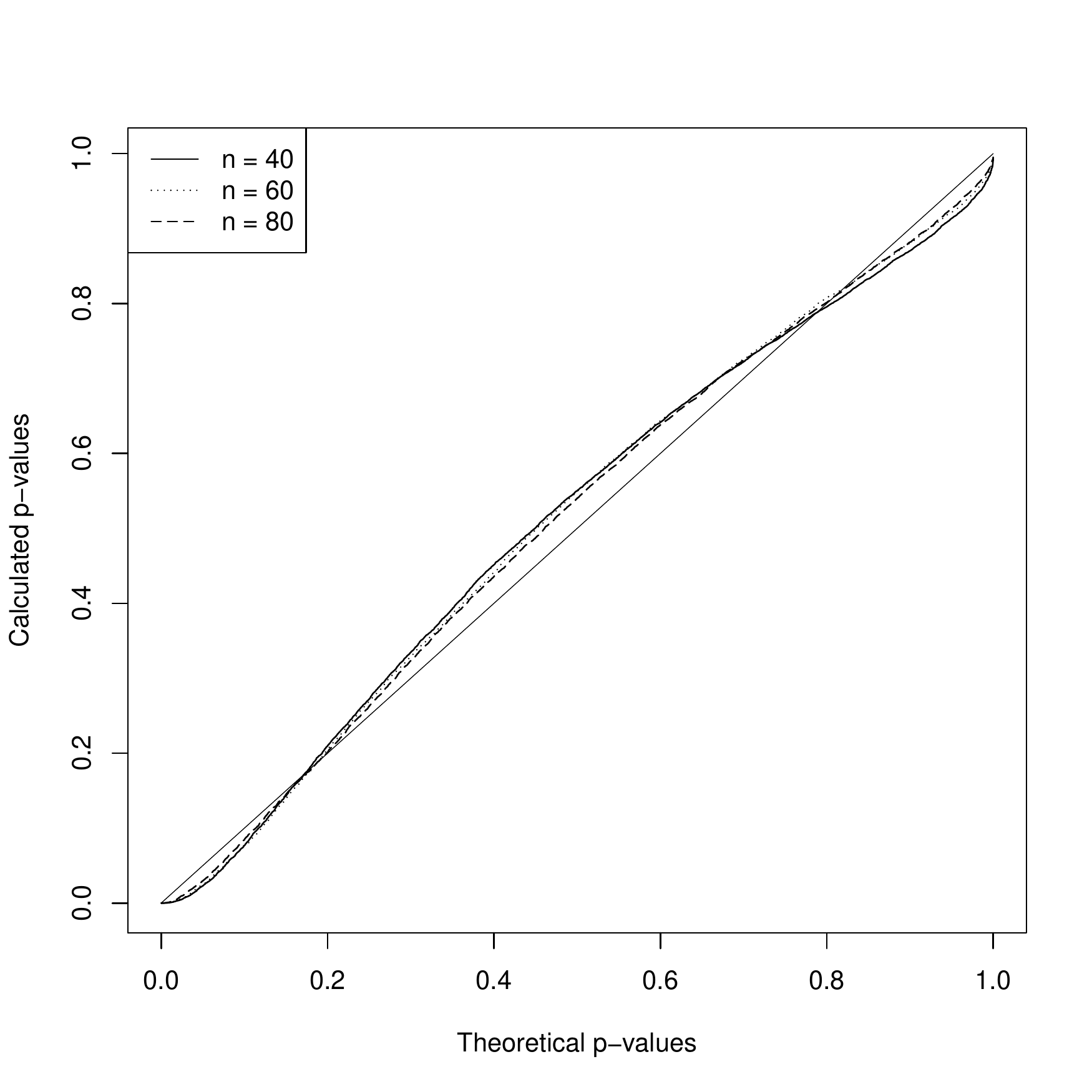} 
 \includegraphics[width=0.3\textwidth]{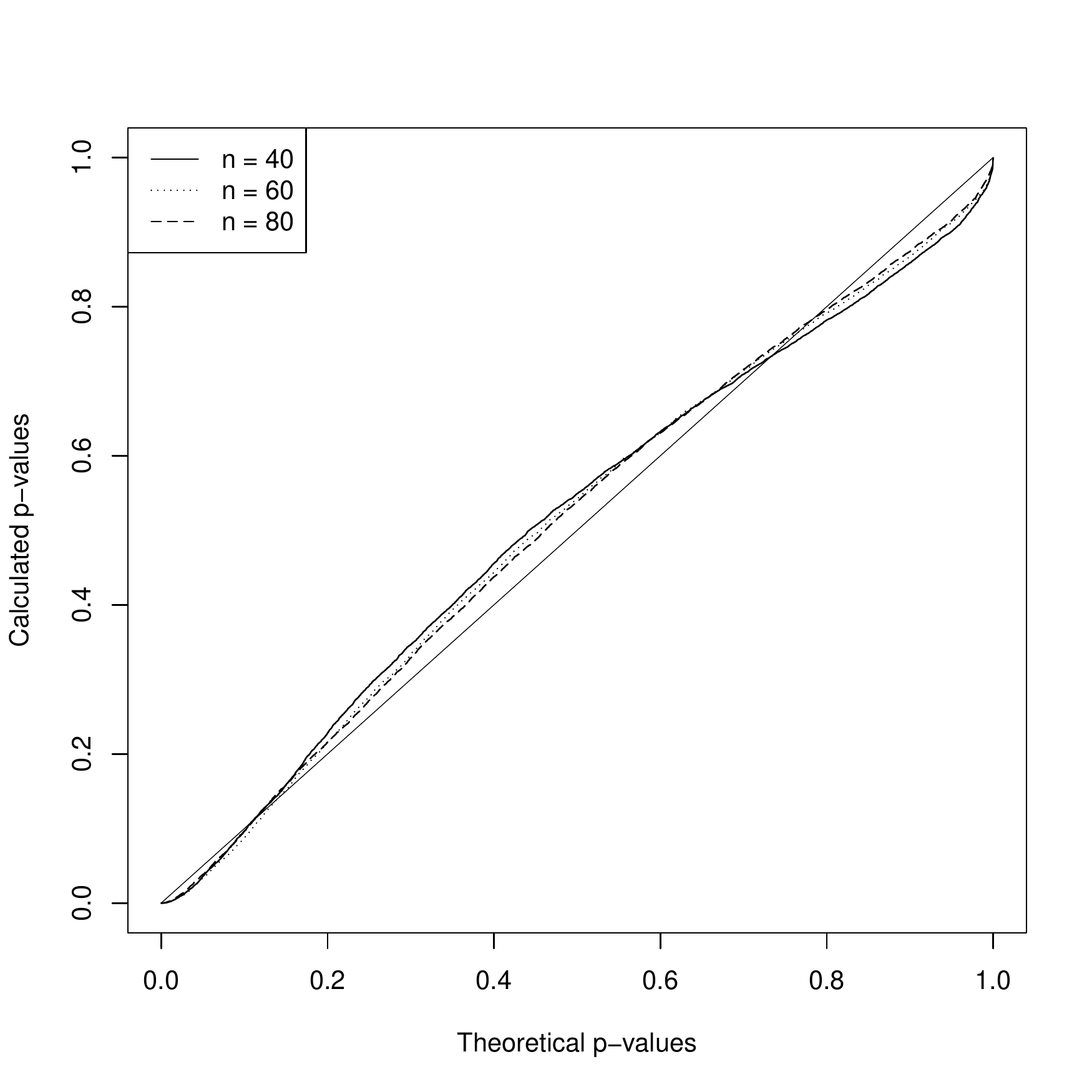} 
 \caption{QQ-plot of the p-values for the two sample case with $M = 1,2$ and $3$ respectively.}
 \label{fig:qqplotmodel12}
\end{figure}

\pagebreak
\bibliographystyle{plain}
\bibliography{references}

\end{document}